\documentclass[a4paper,11pt,reqno]{amsart}

\usepackage[utf8]{inputenc}
\usepackage{mathrsfs}

\usepackage{geometry}
\usepackage[hidelinks]{hyperref}

\usepackage{enumitem,mathtools,mathabx}

\mathtoolsset{showonlyrefs}

\theoremstyle{plain} 
\newtheorem{theorem}{Theorem}[section]
\newtheorem{proposition}[theorem]{Proposition}

\theoremstyle{remark} 
\newtheorem{remark}[theorem]{Remark}

\newcommand{\f}[2]{\frac{#1}{#2}}
\newcommand{\sn}[1]{\mathring{#1}}
\newcommand{\R}{\mathbb{R}}
\newcommand{\C}{\mathbb{C}}
\newcommand{\e}{\mathrm{e}}
\newcommand{\ep}{\varepsilon}
\newcommand{\F}{\mathscr{F}}
\newcommand{\hd}[1]{\mathcal{H}_\alpha^#1}
\newcommand{\Ht}[2]{\mathcal{H}_{\alpha(#2)}^#1}
\newcommand{\hn}[1]{\mathcal{H}_{\beta,\sigma}^#1}
\newcommand{\hS}[1]{\mathcal{H}_{\beta,\sigma}^{#1}}
\newcommand{\hs}[1]{\mathcal{H}_{\beta,\sigma}^{#1,s}}
\newcommand{\en}[1]{E_{\beta,\sigma}^#1}
\newcommand{\En}[1]{\mathcal{E}_{\beta,\sigma}^#1}
\newcommand{\M}[1]{\mathcal{M}^#1}
\newcommand{\qd}[1]{\mathcal{Q}_\alpha^#1}
\newcommand{\Ud}[1]{\mathcal{U}_#1}
\newcommand{\Kd}[1]{\mathcal{K}_#1}
\newcommand{\D}{\mathrm{dom}}
\newcommand{\Ld}[2]{L^#1(\R^#2)}
\newcommand{\Hd}[2]{H^#1(\R^#2)}
\newcommand{\Hdloc}[2]{H_{loc}^#1(\R^#2)}
\newcommand{\Ho}[2]{\sn{H}^#1(\R^#2)}
\newcommand{\Lu}[1]{L^#1(\R)}
\newcommand{\Hu}[1]{H^#1(\R)}
\newcommand{\Gd}[2]{\mathcal{G}_#1^#2}
\newcommand{\GPd}[2]{\mathscr{G}_#1^#2}
\newcommand{\Gs}[2]{\mathcal{G}_#1^{#2,s}}
\newcommand{\td}[3]{\theta_#1^#2\left(#3\right)}
\newcommand{\tn}[3]{T_{#1,#2}^#3}
\newcommand{\un}[3]{u_{#1,#2,\omega}^#3}
\newcommand{\unu}[3]{u_{#1,#2,1}^#3}
\newcommand{\psin}[3]{\psi_{#1,#2,\omega}^#3}
\renewcommand{\S}{\mathbb{S}^2}
\newcommand{\ds}{\,ds}
\newcommand{\dx}{\,dx}
\newcommand{\dy}{\,dy}
\newcommand{\dsy}{\,d\S(y)}

\renewcommand{\i}{\imath}

\begin{document}

\title[A general review on the NLS equation with point-concentrated nonlinearity]{A general review on the NLS equation with point-concentrated nonlinearity}

\author[L. Tentarelli]{Lorenzo Tentarelli}
\address[L. Tentarelli]{Politecnico di Torino, Dipartimento di Scienze Matematiche ``G.L. Lagrange'', Corso Duca degli Abruzzi 24, 10129, Torino, Italy.} 
\email{lorenzo.tentarelli@polito.it}

\date{\today}

\begin{abstract} 
 The paper presents a complete (to the best of the author's knowledge) overview on the existing literature concerning the NLS equation with point-concentrated nonlinearity. Precisely, it mainly covers the following topics: definition of the model, weak and strong local well-posedness, global well-posedness, classification and stability (orbital and asymptotic) of the standing waves, blow-up analysis and derivation from the standard NLS equation with shrinking potentials. Also some related problem is mentioned.
\end{abstract}

\maketitle

\vspace{-.5cm}
{\footnotesize AMS Subject Classification: 35Q40, 35Q55, 35R06, 81Q99}
\smallskip

{\footnotesize Keywords: NLS equation, concentrated nonlinearity, delta potentials, well-posedness, standing waves, stability, blow-up, point-like limit.}
    
%%%%%%%%%%%%%%%%%%%%%%%%%%%%%%%%%%%%%%%%%%%%%%%%%%%%%%%%%%%%%%%%%%%%%%%%%%%%%%%%
%%%%%%%%%%%%%%%%%%%%%%%%%%%%%%%%%%%%%%%%%%%%%%%%%%%%%%%%%%%%%%%%%%%%%%%%%%%%%%%%
%%%%%%%%%%%%%%%%%%%%%%%%%%%%%%%%%%%%%%%%%%%%%%%%%%%%%%%%%%%%%%%%%%%%%%%%%%%%%%%%
%%%%%%%%%%%%%%%%%%%%%%%%%%%%%%%%%%%%%%%%%%%%%%%%%%%%%%%%%%%%%%%%%%%%%%%%%%%%%%%%
%%%%%%%%%%%%%%%%%%%%%%%%%%%%%%%%%%%%%%%%%%%%%%%%%%%%%%%%%%%%%%%%%%%%%%%%%%%%%%%%
%%%%%%%%%%%%%%%%%%%%%%%%%%%%%%%%%%%%%%%%%%%%%%%%%%%%%%%%%%%%%%%%%%%%%%%%%%%%%%%%
%%%%%%%%%%%%%%%%%%%%%%%%%%%%%%%%%%%%%%%%%%%%%%%%%%%%%%%%%%%%%%%%%%%%%%%%%%%%%%%%
%%%%%%%%%%%%%%%%%%%%%%%%%%%%%%%%%%%%%%%%%%%%%%%%%%%%%%%%%%%%%%%%%%%%%%%%%%%%%%%%
%%%%%%%%%%%%%%%%%%%%%%%%%%%%%%%%%%%%%%%%%%%%%%%%%%%%%%%%%%%%%%%%%%%%%%%%%%%%%%%%
%%%%%%%%%%%%%%%%%%%%%%%%%%%%%%%%%%%%%%%%%%%%%%%%%%%%%%%%%%%%%%%%%%%%%%%%%%%%%%%%
%%%%%%%%%%%%%%%%%%%%%%%%%%%%%%%%%%%%%%%%%%%%%%%%%%%%%%%%%%%%%%%%%%%%%%%%%%%%%%%%
%%%%%%%%%%%%%%%%%%%%%%%%%%%%%%%%%%%%%%%%%%%%%%%%%%%%%%%%%%%%%%%%%%%%%%%%%%%%%%%%

\section{Introduction}
\label{sec-intro}

The standard NonLinear Schr\"odinger (NLS) equation, i.e.,
\begin{equation}
\label{eq-NLS}
 \i\f{\partial\psi}{\partial t}=-\Delta\psi+\beta|\psi|^{2\sigma}\psi,\qquad \beta\in\R\setminus\{0\},\quad\sigma>0,\qquad\text{on}\quad\R^+\times\R^d,
\end{equation}
is known to play a relevant role in several sectors of physics, where it appears as an effective evolution equation describing the evolution of a microscopic system on a macroscopic or mesoscopic scale. A typical example is provided by Bose-Einstein Condensates (BEC), whose behavior is ``well approximated'' in some suitable sense by the solutions of \eqref{eq-NLS} (see, e.g., \cite{DGPS-99}). There are, however, other examples in which the physical meaning of the NLS equation is totally different, such as, for instance, the propagation of light in nonlinear optics, the behavior of water or plasma waves, the signal transmission through neurons (FitzHugh-Nagumo model), etc. (see \cite{M-05} and references therein).

This paper is, to the best of our knowledge, a complete overview on the literature concerning one of the possible singular perturbations of \eqref{eq-NLS}: the NLS equation with \emph{point-concentrated} nonlinearity (or Schr\"odinger equation with \emph{nonlinear delta potential}). Such equation reads as \eqref{eq-NLS}, but with the nonlinearity formally multiplied times a \emph{Dirac delta measure} based at the origin, i.e.
\begin{equation}
\label{eq-formalDNLS}
 \i\f{\partial\psi}{\partial t}=-\Delta\psi+\beta|\psi|^{2\sigma}\psi\delta_0,\qquad \beta\in\R\setminus\{0\},\quad\sigma>0,\qquad\text{on}\quad\R^+\times\R^d.
\end{equation}

The idea of introducing point perturbations of delta-type in Quantum Mechanics traces back to 1936 and is due to E. Fermi, though in a completely different context as he proposed linear point perturbations of delta-type of the Schr\"odinger equation (i.e., $\sigma=0$ above) to model the interaction between a slow neutron and a fixed atom (see \cite{F-36}).

The interest for the nonlinear model depicted in \eqref{eq-formalDNLS}, on the contrary, is much more recent (and related to the applications of the standard NLS equation). It has been introduced in the 1990s to describe phenomena mostly related to solid state and condensed matter physics: charge accumulation in semiconductor interfaces or heterostructures (\cite{BKB-96,JPS-95,MA-93,MB-02,N-98,PJC-91}), nonlinear propagation in a Kerr-type media in presence of localized defects (\cite{SKB-99,SKBRC-01,Y-05}), BEC in optical lattices where an isolated defect is generated by a focused laser beam (\cite{DM-11,LKMF-11}). Other applications are also suggested in acoustic, conventional and high-$T_c$ superconductivity, light propagation in photonic crystals, etc. (see, e.g., \cite{SKBRC-01} and references therein). However, a rigorous mathematical derivation of the model in these contexts is still missing. 

Still from the perspective of modeling BEC in the presence of defects or impurities (whose spatial scale is supposed to be much smaller than the dispersion of the wave function) it is worth mentioning that also another equation involving delta potentials has been proposed in recent years (\cite{SM-20,SCMS-20}):
\begin{equation}
\label{eq-altraDNLS}
 \i\f{\partial\psi}{\partial t}=(-\Delta+\alpha\delta_0)\psi+\beta|\psi|^{2\sigma}\psi,\qquad \alpha\in\R,\quad\beta\in\R\setminus\{0\},\quad\sigma>0,\qquad\text{on}\quad\R^+\times\R^d.
\end{equation}
Results on \eqref{eq-altraDNLS}, which are not presented in this review, are discussed by \cite{AN-09,AN-13,ANV-13,DH-09,FJ-08,FOO-08,GHW-04,HMZ-07,LFFKS-08,MN-21} in dimension one and by \cite{ABCT-22a,ABCT-22b,CFN-21,CFN-23,FGI-22,FN-22} in dimensions two and three (while \cite{FT-00} addresses the circle, \cite{BC-22} addresses the half-line and \cite{ABD-22,BD-21,BD-22} concern some first studies on a mixed model between \eqref{eq-formalDNLS} and \eqref{eq-altraDNLS}).

Finally, we mention that models involving delta potentials, as the one addressed by this paper, are rigorously defined only for $d=1,\,2,\,3$ (more details on this point will be provided in Section \ref{sec-model}). Hence, we will tacitly assume this restriction on the space dimension $d$ of $\R^d$ throughout the review.

\begin{remark}
In this paper we only focus on the case where a \emph{single} nonlinear delta potential is based at the \emph{origin}, i.e. $\beta|\psi|^{2\sigma}\psi\delta_0$. However, anything can be suitably adapted to the case of a delta based in other points of $\R^d$. Moreover, some results can be adapted to the case of finitely many deltas. We decided to limit ourselves to \eqref{eq-formalDNLS} for the sake of simplicity, but, nevertheless, we will mention the cases in which the results are known to extend to the more general frameworks.
\end{remark}

%%%%%%%%%%%%%%%%%%%%%%%%%%%%%%%%%%%%%%%%%%%%%%%%%%%%%%%%%%%%%%%%%%%%%%%%%%%%%%%%
%%%%%%%%%%%%%%%%%%%%%%%%%%%%%%%%%%%%%%%%%%%%%%%%%%%%%%%%%%%%%%%%%%%%%%%%%%%%%%%%
%%%%%%%%%%%%%%%%%%%%%%%%%%%%%%%%%%%%%%%%%%%%%%%%%%%%%%%%%%%%%%%%%%%%%%%%%%%%%%%%
%%%%%%%%%%%%%%%%%%%%%%%%%%%%%%%%%%%%%%%%%%%%%%%%%%%%%%%%%%%%%%%%%%%%%%%%%%%%%%%%
%%%%%%%%%%%%%%%%%%%%%%%%%%%%%%%%%%%%%%%%%%%%%%%%%%%%%%%%%%%%%%%%%%%%%%%%%%%%%%%%
%%%%%%%%%%%%%%%%%%%%%%%%%%%%%%%%%%%%%%%%%%%%%%%%%%%%%%%%%%%%%%%%%%%%%%%%%%%%%%%%

\subsection{Organization of the paper}
\label{subsec:org}

The paper is organized as follows.
\begin{itemize}
 \item Section \ref{sec-model} addresses the rigorous definition of the model, starting from the linear case (Section \ref{subsec-linear}) to the nonlinear case (Section \ref{subsec-nonlinear}).
 \item Section \ref{sec-LWP} addresses local well-posedness of the Cauchy problem associated with \eqref{eq-formalDNLS}, both in the \emph{weak} and in the \emph{strong} sense.
 \item Section \ref{sec-GWP} addresses global well-posedness of the Cauchy problem associated with \eqref{eq-formalDNLS}.
 \item Section \ref{sec-standing} addresses the classification and the stability (\emph{orbital} in Section \ref{subsec-orbital}, \emph{asymptotic} in Section \ref{subsec-asymptotic}) of the standing waves of \eqref{eq-formalDNLS}.
 \item Section \ref{sec-blow-up} addresses blow-up analysis, from the first results (Section \ref{subsec-seminal}) to the the more refined ones (Section \ref{subsec-further}), also mentioning the question of pseudoconformal invariance (Section \ref{subsec-pseudo}).
 \item Section \ref{sec-derivation} addressed the derivation of \eqref{eq-formalDNLS} from the standard NLS equation in dimension one (Section \ref{subsec-1d}) and three (Section \ref{subsec-3d}).
 \item Section \ref{sec-connected} briefly mentions some connected problems.
\end{itemize}

%%%%%%%%%%%%%%%%%%%%%%%%%%%%%%%%%%%%%%%%%%%%%%%%%%%%%%%%%%%%%%%%%%%%%%%%%%%%%%%%
%%%%%%%%%%%%%%%%%%%%%%%%%%%%%%%%%%%%%%%%%%%%%%%%%%%%%%%%%%%%%%%%%%%%%%%%%%%%%%%%
%%%%%%%%%%%%%%%%%%%%%%%%%%%%%%%%%%%%%%%%%%%%%%%%%%%%%%%%%%%%%%%%%%%%%%%%%%%%%%%%
%%%%%%%%%%%%%%%%%%%%%%%%%%%%%%%%%%%%%%%%%%%%%%%%%%%%%%%%%%%%%%%%%%%%%%%%%%%%%%%%
%%%%%%%%%%%%%%%%%%%%%%%%%%%%%%%%%%%%%%%%%%%%%%%%%%%%%%%%%%%%%%%%%%%%%%%%%%%%%%%%
%%%%%%%%%%%%%%%%%%%%%%%%%%%%%%%%%%%%%%%%%%%%%%%%%%%%%%%%%%%%%%%%%%%%%%%%%%%%%%%%

\subsection*{Fundings
%and acknowledgements
} The author has been partially supported by the INdAM GNAMPA project 2022 ``Modelli matematici con singolarit\`a per fenomeni di interazione'' (CUP E55F2200027\-0001).

% \subsection*{Data availability statement}
% Data sharing not applicable to this article as no datasets were generated or analysed during the current study.
%  
% \subsection*{Conflict of interest}
% On behalf of all authors, the corresponding author states that there is no conflict of interest.

%%%%%%%%%%%%%%%%%%%%%%%%%%%%%%%%%%%%%%%%%%%%%%%%%%%%%%%%%%%%%%%%%%%%%%%%%%%%%%%%
%%%%%%%%%%%%%%%%%%%%%%%%%%%%%%%%%%%%%%%%%%%%%%%%%%%%%%%%%%%%%%%%%%%%%%%%%%%%%%%%
%%%%%%%%%%%%%%%%%%%%%%%%%%%%%%%%%%%%%%%%%%%%%%%%%%%%%%%%%%%%%%%%%%%%%%%%%%%%%%%%
%%%%%%%%%%%%%%%%%%%%%%%%%%%%%%%%%%%%%%%%%%%%%%%%%%%%%%%%%%%%%%%%%%%%%%%%%%%%%%%%
%%%%%%%%%%%%%%%%%%%%%%%%%%%%%%%%%%%%%%%%%%%%%%%%%%%%%%%%%%%%%%%%%%%%%%%%%%%%%%%%
%%%%%%%%%%%%%%%%%%%%%%%%%%%%%%%%%%%%%%%%%%%%%%%%%%%%%%%%%%%%%%%%%%%%%%%%%%%%%%%%
%%%%%%%%%%%%%%%%%%%%%%%%%%%%%%%%%%%%%%%%%%%%%%%%%%%%%%%%%%%%%%%%%%%%%%%%%%%%%%%%
%%%%%%%%%%%%%%%%%%%%%%%%%%%%%%%%%%%%%%%%%%%%%%%%%%%%%%%%%%%%%%%%%%%%%%%%%%%%%%%%
%%%%%%%%%%%%%%%%%%%%%%%%%%%%%%%%%%%%%%%%%%%%%%%%%%%%%%%%%%%%%%%%%%%%%%%%%%%%%%%%
%%%%%%%%%%%%%%%%%%%%%%%%%%%%%%%%%%%%%%%%%%%%%%%%%%%%%%%%%%%%%%%%%%%%%%%%%%%%%%%%
%%%%%%%%%%%%%%%%%%%%%%%%%%%%%%%%%%%%%%%%%%%%%%%%%%%%%%%%%%%%%%%%%%%%%%%%%%%%%%%%
%%%%%%%%%%%%%%%%%%%%%%%%%%%%%%%%%%%%%%%%%%%%%%%%%%%%%%%%%%%%%%%%%%%%%%%%%%%%%%%%

\section{Definition of the model}
\label{sec-model}

In order to give a rigorous meaning to \eqref{eq-formalDNLS} it is necessary to start from the linear analogous ($\sigma=0$) and, then, go back to the nonlinear problem. According to what we said before, we limit ourselves to present the setting in the case of a single delta based at the origin.

%%%%%%%%%%%%%%%%%%%%%%%%%%%%%%%%%%%%%%%%%%%%%%%%%%%%%%%%%%%%%%%%%%%%%%%%%%%%%%%%
%%%%%%%%%%%%%%%%%%%%%%%%%%%%%%%%%%%%%%%%%%%%%%%%%%%%%%%%%%%%%%%%%%%%%%%%%%%%%%%%
%%%%%%%%%%%%%%%%%%%%%%%%%%%%%%%%%%%%%%%%%%%%%%%%%%%%%%%%%%%%%%%%%%%%%%%%%%%%%%%%
%%%%%%%%%%%%%%%%%%%%%%%%%%%%%%%%%%%%%%%%%%%%%%%%%%%%%%%%%%%%%%%%%%%%%%%%%%%%%%%%
%%%%%%%%%%%%%%%%%%%%%%%%%%%%%%%%%%%%%%%%%%%%%%%%%%%%%%%%%%%%%%%%%%%%%%%%%%%%%%%%
%%%%%%%%%%%%%%%%%%%%%%%%%%%%%%%%%%%%%%%%%%%%%%%%%%%%%%%%%%%%%%%%%%%%%%%%%%%%%%%%

\subsection{The linear model}
\label{subsec-linear}

The most suitable way to give a rigorous meaning to a perturbation of delta type of the Laplacian, namely to an operator of the type $-\Delta+\alpha\delta_0$, with $\alpha\in\R$, is applying the theory of \emph{self-adjoint extensions} of symmetric operators to $-\Delta_{|C_0^\infty(\R^d\setminus\{0\})}$ (see, e.g., \cite{AGH-KH-88} for a general dissertation and \cite{A-00} for some basics). In particular, one can prove that nontrivial extensions exist only for $d=1,\,2,\,3$ and that, for any fixed $\alpha\in\R$, the unique self-adjoint realization of $-\Delta+\alpha\delta_0$ is given by the operator $\hd{d}:\Ld{2}{d}\to\Ld{2}{d}$ with \emph{domain}
\begin{multline}
\label{eq-dom}
 \D(\hd{d}):=\left\{u\in\Ld{2}{d}:\exists\lambda>0,\,\exists q\in\C\:\text{ such that }\right.\\
 \left.u-\kappa_d(\alpha)q\Gd{\lambda}{d}=:\phi_\lambda\in\Hd{2}{d}\:\text{ and }\:\phi_\lambda(0)=\td{\lambda}{d}{\alpha}q\right\}
\end{multline}
and \emph{action}
\begin{equation}
 \label{eq-act}
 \hd{d}u:=-\Delta\phi_\lambda-\kappa_d(\alpha)q\lambda\Gd{\lambda}{d},\qquad\forall u\in\D(\hd{d}),
\end{equation}
where
\begin{equation}
 \label{eq-theta}
 \kappa_d(\alpha):=
 \left\{
 \begin{array}{ll}
 -\alpha, & \text{if}\quad d=1,\\[.4cm]
 1, & \text{if}\quad d=2,\,3,
 \end{array}
 \right. 
 \qquad
 \td{\lambda}{d}{\alpha}:=
 \left\{
 \begin{array}{ll}
  \displaystyle \f{2\sqrt{\lambda}}{\alpha+2\sqrt{\lambda}}, & \text{if}\quad d=1,\\[.4cm]
  \displaystyle \alpha+\f{\log\left(\f{\sqrt{\lambda}}{2}\right)+\gamma}{2\pi}, & \text{if}\quad d=2,\\[.4cm]
  \displaystyle \alpha+\f{\sqrt{\lambda}}{4\pi}, & \text{if}\quad d=3,
 \end{array}
 \right.
\end{equation}
with $\gamma$ the \emph{Euler-Mascheroni} constant, and $\Gd{\lambda}{d}$ is the \emph{Green's function} of $-\Delta+\lambda$, i.e.
\begin{equation}
 \label{eq-green}
 \Gd{\lambda}{d}(x)=\F^{-1}\bigg(\f{1}{(2\pi)^{d/2}(|k|^2+\lambda)}\bigg)(x)=
 \left\{
 \begin{array}{ll}
  \displaystyle \f{\e^{-\sqrt{\lambda}|x|}}{2\sqrt{\lambda}}, & \text{if}\quad d=1,\\[.4cm]
  \displaystyle \f{K_0(\sqrt{\lambda}|x|)}{2\pi}, & \text{if}\quad d=2,\\[.4cm]
  \displaystyle \f{\e^{-\sqrt{\lambda}|x|}}{4\pi|x|}, & \text{if}\quad d=3,
 \end{array}
 \right.
\end{equation}
with $\F$ denoting the unitary Fourier transform of $\R^d$ (i.e., $\F(h)[k]:=(2\pi)^{-d/2}\int_{\R^d}\e^{-\i k\cdot x}h(x)\dx$, whenever $h\in \Ld{1}{d}$) and $K_0$ denoting the modified Bessel function of the second kind of order 0, also known as \emph{Macdonald function} (see, e.g., \cite[Section 9.6]{AS-65}). Note that $\D(\hd{d})$ is a Hilbert space if endowed with the \emph{graph norm} $\|\cdot\|_{\D(\hd{d})}^2:=\|\cdot\|_{\Ld{2}{d}}^2+\|\hd{d}\,\cdot\|_{\Ld{2}{d}}^2$.

\begin{remark}
 Note that, when $d=1$, $\alpha=0$ represents the free Laplacian. On the contrary, it can be recovered in $d=2,\,3$ only letting $\alpha\to\infty$.
\end{remark}

As a consequence of \eqref{eq-dom}, any function in $\D(\hd{d})$ admits a decomposition in  a \emph{regular part} $\phi_\lambda$, on which the operator acts as the standard Laplacian, and a \emph{singular part} $q\Gd{\lambda}{d}$, on which the operator acts as the multiplication times $-\lambda$. In addition, we recall that the two components are bound by the so-called \emph{boundary condition} $\phi_{\lambda}(0)=\td{\lambda}{d}{\alpha}q$ and that the strength $q$ of the singular part is usually called \emph{charge}. The charge is uniquely determined for any $u\in\D(\hd{d})$ and different charges identify different functions $u\in\D(\hd{d})$ (see, e.g., \cite[Remark 2.1]{ABCT-22a} for the case $d=2$). On the contrary, $\lambda$ is a dumb parameter in the sense that every function of $\D(\hd{d})$ admits an equivalent decomposition for any positive value of $\lambda$, due to the fact that the difference between Green's functions with different values of $\lambda$ is in $\Hd{2}{d}$.

As a further evidence that the parameter $\lambda$ cannot actually affect the definition of $\hd{d}$, it is possible to find formulations equivalent to \eqref{eq-dom} and \eqref{eq-act} that do not involve $\lambda$. In other words, one can check that $\D(\hd{d})$ is equal to
 \begin{gather}
  \label{eq-1d0} \displaystyle \left\{u\in\Hu{1}\cap H^2(\R\setminus\{0\}):u'(0^+)-u'(0^-)=\alpha u(0)\right\},\quad \text{if}\quad d=1,\\[.4cm]
  \label{eq-2d0} \displaystyle \left\{u\in\Ld{2}{2}:\exists q\in\C\:\text{ s.t. }u-q\Gd{0}{2}=:\phi\in\Hdloc{2}{2}\cap\Ho{2}{2},\:\phi(0)=\alpha q\right\},\quad \text{if}\quad d=2,\\[.4cm]
  \label{eq-3d0} \displaystyle \left\{u\in\Ld{2}{3}:\exists q\in\C\:\text{ s.t. }u-q\Gd{0}{3}=:\phi\in\Hdloc{2}{3}\cap\Ho{1}{3}\cap\Ho{2}{3},\:\phi(0)=\alpha q\right\},\quad \text{if}\quad d=3,
\end{gather}
where
\begin{equation}
 \label{eq-greendue}
 \Gd{0}{d}(x)=\F^{-1}\bigg(\f{1}{(2\pi)^{d/2}|k|^2}\bigg)(x)=
 \left\{
 \begin{array}{ll}
  \displaystyle -\f{\log|x|}{2\pi}, & \text{if}\quad d=2,\\[.4cm]
  \displaystyle \f{1}{4\pi|x|}, & \text{if}\quad d=3,
 \end{array}
 \right.
\end{equation}
and that
\begin{equation}
 \label{eq-actdue}
 \hd{d}u=
 \left\{
 \begin{array}{lll}
  \displaystyle -\f{d^2u}{dx^2}, & \text{in}\quad \R\setminus\{0\}, & \text{if}\quad d=1,\\[.4cm]
  \displaystyle -\Delta \phi, & \text{in}\quad \R^d\setminus\{0\}, & \text{if}\quad d=2,\,3,
 \end{array}
 \right.
 \qquad\forall u\in\D(\hd{d}).
\end{equation}

In view of the previous definitions, it is not evident that $\hd{d}$ is the proper (self-adjoint) realization of the delta perturbation of the Laplacian. However, easy computations yields that, in the sense of distributions,
\[
 \hd{d}u=-\Delta u-\kappa_d(\alpha)q\delta_0,\qquad\forall u\in\D(\hd{d}).
\]
Note also that, when $d=1$, there results $q=u(0)$, so that one can rigorously claim that $\hd{d}=-\f{d^2}{dx^2}+\alpha\delta_0$. Thus, the construction provided before is not strictly necessary for $d=1$. This is due to the fact that $-\f{d^2}{dx^2}+\alpha\delta_0$ is a form bounded perturbation of the Laplacian, which makes the one-dimensional case considerably different from the other ones because of the KLMN Theorem (see, e.g., \cite[Theorem X.17]{RS-75}).

\begin{remark}
It is worth mentioning that the one-dimensional case is special also from another point of view. Indeed, only in this case, the one described above is not the unique class of nontrivial self-adjoint extension of $-\Delta_{|C_0^\infty(\R\setminus\{0\})}$. In fact, when $d=1$ there is a \emph{four parameters} family of self-adjoint extensions (while in $d=2,\,3$ there is only a \emph{one parameter} family of self-adjoint extensions). However, these further generalizations are not of delta type and go beyond the aims of the present review, so that are omitted.
\end{remark}

\begin{remark}
 Another way to obtain delta perturbations of the Laplacian is by studying the asymptotic behavior of perturbations of the Laplacian consisting of smooth suitably shrinking potentials (both local and nonlocal). This can be found, for instance, in \cite{AGH-KH-88} and brings exactly to the operator $\hd{d}$ introduced before. 
\end{remark}

The \emph{spectrum} $\sigma(\hd{d})$ varies with the dimension as follows
\begin{equation}
 \label{eq-sp}
 \sigma(\hd{d}):=
 \left\{
 \begin{array}{cl}
  \displaystyle
  \left\{
  \begin{array}{ll}
  \displaystyle [0,+\infty), & \text{if}\quad \alpha\geq0,\\[.4cm]
  \displaystyle \{\ell_\alpha^1\}\cup[0,+\infty), & \text{if}\quad \alpha<0,
  \end{array}
  \right\}
  & \text{if}\quad d=1,\\[.8cm]
  \displaystyle \{\ell_\alpha^2\}\cup[0,+\infty), & \text{if}\quad d=2,\\[.4cm]
  \left\{
  \begin{array}{ll}
  \displaystyle [0,+\infty), & \text{if}\quad \alpha\geq0,\\[.4cm]
  \displaystyle \{\ell_\alpha^3\}\cup[0,+\infty), & \text{if}\quad \alpha<0,
  \end{array}
  \right\}
  & \text{if}\quad d=3,
 \end{array}
 \right.
\end{equation}
with
\begin{equation}
 \label{eq-eigen}
 \ell_\alpha^d:=\left\{
 \begin{array}{ll}
  \displaystyle -\alpha^2/4, & \text{if}\quad d=1,\quad\text{and}\quad\alpha<0,\\[.4cm]
  \displaystyle -4\e^{-4\pi\alpha-2\gamma} & \text{if}\quad d=2,\quad\text{and}\quad\alpha\in\R,\\[.4cm]
  \displaystyle -16\pi^2\alpha^2, & \text{if}\quad d=3,\quad\text{and}\quad\alpha<0.
 \end{array}
 \right.
\end{equation}
where, as usual in the context of the operator theory, the presence or the absence of the negative \emph{eigenvalue} $\ell_\alpha^d$ distinguish between an \emph{attractive} delta potential and a \emph{repulsive} delta potential (note that in $d=2$ the delta potential is always attractive). In the cases in which $\ell_\alpha^d$ is present its \emph{eigenspace} is spanned by $\Gd{{-\ell_\alpha^d}}{d}$.

Moreover, the \emph{quadratic form} associated with $\hd{d}$ reads
\begin{equation}
\label{eq-qform}
 \qd{d}(u):=\left\{
 \begin{array}{ll}
  \displaystyle \left\|\f{du}{dx}\right\|_{\Lu{2}}^2+\alpha|q|^2, & \text{if}\quad d=1\\[.4cm]
  \displaystyle \|\nabla\phi_\lambda\|_{\Ld{2}{d}}^2+\lambda(\|\phi_\lambda\|_{\Ld{2}{d}}^2-\|u\|_{\Ld{2}{d}}^2)+\td{\lambda}{d}{\alpha}|q|^2, & \text{if}\quad d=2,\,3,
 \end{array}
 \right.
\end{equation}
for every $u\in\D(\qd{d})$, where
\begin{equation}
\label{eq-formdom}
 \D(\qd{d})=V_d:=\left\{
 \begin{array}{ll}
  \displaystyle \Hu{1}, & \text{if}\quad d=1,\\[.4cm]
  \displaystyle \{u\in\Ld{2}{d}:u-q\Gd{\lambda}{d}=:\phi_\lambda\in\Hd{1}{d}\}, & \text{if}\quad d=2,\,3.
 \end{array}
 \right.
\end{equation}
We introduced the notation $V_d$ for the \emph{form domain} to underline that, in contrast to the operator domain, it does not depend on $\alpha$. In addition, $V_d\supsetneq \Hd{1}{d}$, when $d\neq1$, and it is a Hilbert space when endowed with the norm
\[
 \|u\|_{V_d}^2:=\left\{
 \begin{array}{ll}
  \displaystyle \|u\|_{\Hu{1}}^2, & \text{if}\quad d=1,\\[.4cm]
  \displaystyle \|\phi_\lambda\|_{\Hd{1}{d}}^2+|q|^2\|\Gd{\lambda}{d}\|_{\Ld{2}{d}}^2, & \text{if}\quad d=2,\,3
 \end{array}
 \right.
\]
(equivalent for every fixed $\lambda>0$, when $d=2,3$). In fact, also the three-dimensional case admits an analogous representation for the quadratic form, which arises letting $\lambda\to0$, i.e.
\[
 \qd{3}(u)=\|\nabla\phi\|_{\Ld{2}{3}}^2+\alpha|q|^2,\qquad\forall u\in V_3,
\]
with
\[
 V_3=\{u\in\Ld{2}{3}:u-q\Gd{0}{3}=:\phi\in\Hdloc{1}{3}\cap\Ho{1}{3}\}.
\]
In this case the norm can be equivalently written as $\|u\|_{V_3}^2=\|\nabla\phi\|_{\Ld{2}{3}}^2+|q|^2$. On the contrary, in the two-dimensional case it is not known whether this kind of representation is available due to the low regularity of the regular part when $\lambda=0$. Such a difference is a consequence of the fact that $\Gd{\lambda}{3}\to\Gd{0}{3}$, while $\Gd{\lambda}{2}\to+\infty$, pointwise as $\lambda\to0$.

\begin{remark}
 We also mention that, whenever $\ell_\alpha^d$ exists, it satisfies
 \[
  \ell_\alpha^d=\qd{d}\left(f_d(-\ell_\alpha^d)\,\Gd{{-\ell_\alpha^d}}{d}\right)=\inf_{\substack{u\in V_d\\\|u\|_{\Ld{2}{d}}=1}}\qd{d}(u),
 \]
 with $f_d:\R^+\to\R^+$ given by
 \[
  f_d(\lambda):=\left\{
 \begin{array}{ll}
  \displaystyle 2\sqrt[4]{\lambda^3}, & \text{if}\quad d=1,\\[.4cm]
  \displaystyle 2\sqrt{\pi\lambda} & \text{if}\quad d=2,\\[.4cm]
  \displaystyle 2\sqrt[4]{4\pi^2\lambda}, & \text{if}\quad d=3,.
 \end{array}
 \right.
\]
\end{remark}

Finally, we recall that, since $\hd{d}$ is self-adjoint, as a consequence of Stone's theorem, for every $\psi_0\in\D(\hd{d})$, there exists a unique function $\psi\in C_{loc}^0([0,+\infty);\D(\hd{d}))\cap C_{loc}^1([0,+\infty);\Ld{2}{d})$ which \emph{strongly} solves
\begin{equation}
\label{eq-cau-lin-st}
 \left\{
 \begin{array}{l}
 \displaystyle \i\f{\partial\psi}{\partial t}=\hd{d}\psi\\[.4cm]
 \displaystyle \psi(0,\cdot)=\psi_0
 \end{array}
 \right.
\end{equation}
in the sense that the former is satisfied as an equality in $C_{loc}^0([0,+\infty);\Ld{2}{d})$ and the latter as an equality in $\D(\hd{d})$. On the other hand, as the first equation in \eqref{eq-cau-lin-st} is invariant under \emph{gauge transformations} and \emph{time translations}, the associated \emph{mass}, i.e.
\begin{equation}
 \label{eq-mass}
 M(t)=M(\psi(t,\cdot)):=\|\psi(t,\cdot)\|_{\Ld{2}{d}}^2,
\end{equation}
and \emph{energy}, i.e. $\f{1}{2}\qd{d}(\psi(t,\cdot))$ in the linear case, are preserved along the flow. Moreover, there exists a $\Ld{2}{d}$-strongly continuous unitary group $(U_d(t))_{t}$ such that $\psi(t,\cdot)=U_d(t)\psi_0$. However, since the explicit form of the kernel of $U_d(t)$ is not easily manageable (see \cite{AGH-KH-88}), it is customary to represent the solution as
\begin{equation}
 \label{eq-lin-duhamel}
 \psi(t,x)=(\Ud{d}(t)\psi_0)(x)+\i\int_0^t\Ud{d}(t-s,x)\kappa_d(\alpha)q(s)\ds,
\end{equation}
where $\Ud{d}(t)$ is the \emph{free propagator} of $\R^d$, with kernel
\begin{equation}
 \label{eq-prop}
 \Ud{d}(t,x):=\F^{-1}(\e^{-\i|k|^2t})[x]=\f{\e^{-\f{|x|^2}{4\i t}}}{(4\i\pi t)^{d/2}},
\end{equation}
and $q(\cdot)$ is the unique solution of
\begin{equation}
 \label{eq-lincharge}
 q(t)+\int_0^t\Kd{d}(t-s)c_d(\alpha)q(s)\ds=m_df_d(t),
\end{equation}
with
\begin{equation}
 \label{eq-altripar}
 c_d(\alpha):=
 \left\{
 \begin{array}{ll}
 \displaystyle \f{\e^{\i\pi/4}}{2\sqrt{\pi}}\alpha, & \text{if}\quad d=1,\\[.4cm]
 \displaystyle 4\pi(\td{1}{2}{\alpha}-\i/8), & \text{if}\quad d=2,\\[.4cm]
 \displaystyle 4\sqrt{\pi}\e^{\i\pi/4}\alpha, & \text{if}\quad d=3,
 \end{array}
 \right. 
 \qquad
 m_d:=
 \left\{
 \begin{array}{ll}
  \displaystyle 1, & \text{if}\quad d=1,\\[.4cm]
  \displaystyle 4\pi, & \text{if}\quad d=2,\\[.4cm]
  \displaystyle c_d(1), & \text{if}\quad d=3,
 \end{array}
 \right.
\end{equation}
\medskip
\begin{equation}
 \label{eq-forcing}
 f_d(t):=\left\{
 \begin{array}{ll}
 \displaystyle (\Ud{1}(t)\psi_0)(0), & \text{if}\quad d=1\\[.4cm]
  \displaystyle \int_0^t\Kd{d}(t-s)(\Ud{d}(s)\psi_0)(0)\ds, & \text{if}\quad d=2,\,3,
  \end{array}
  \right.
\end{equation}
and
\begin{equation}
 \label{eq-kernel}
 \Kd{d}(t):=\left\{
 \begin{array}{ll}
 \displaystyle \f{1}{\sqrt{t}}, & \text{if}\quad d=1,\,3,\\[.4cm]
  \displaystyle \int_0^{+\infty}\f{t^{s-1}}{\Gamma(s)}\ds, & \text{if}\quad d=2.
  \end{array}
  \right.
\end{equation}
Since at any time $t$ the solution $q(t)$ of \eqref{eq-lincharge} is the coefficient of the singular part of $\psi(t,\cdot)$, with a little abuse of notation, it is called \emph{charge} and thus \eqref{eq-lincharge} is called \emph{linear charge equation}.

%%%%%%%%%%%%%%%%%%%%%%%%%%%%%%%%%%%%%%%%%%%%%%%%%%%%%%%%%%%%%%%%%%%%%%%%%%%%%%%%
%%%%%%%%%%%%%%%%%%%%%%%%%%%%%%%%%%%%%%%%%%%%%%%%%%%%%%%%%%%%%%%%%%%%%%%%%%%%%%%%
%%%%%%%%%%%%%%%%%%%%%%%%%%%%%%%%%%%%%%%%%%%%%%%%%%%%%%%%%%%%%%%%%%%%%%%%%%%%%%%%
%%%%%%%%%%%%%%%%%%%%%%%%%%%%%%%%%%%%%%%%%%%%%%%%%%%%%%%%%%%%%%%%%%%%%%%%%%%%%%%%
%%%%%%%%%%%%%%%%%%%%%%%%%%%%%%%%%%%%%%%%%%%%%%%%%%%%%%%%%%%%%%%%%%%%%%%%%%%%%%%%
%%%%%%%%%%%%%%%%%%%%%%%%%%%%%%%%%%%%%%%%%%%%%%%%%%%%%%%%%%%%%%%%%%%%%%%%%%%%%%%%

\subsection{From linear to nonlinear}
\label{subsec-nonlinear}

In view of the previous section, the nonlinear problem arises simply letting the ``strength'' $\alpha$ of the delta perturbation depend on the unknown $\psi$, and more precisely on the charge $q$, as follows:
\begin{equation}
 \label{eq-nonlinearity}
 \alpha=\alpha(q):=\beta|q|^{2\sigma}\qquad \beta\in\R\setminus\{0\},\quad\sigma>0.
\end{equation}
As a consequence, the associated Cauchy problem reads
\begin{equation}
\label{eq-cau-st}
 \left\{
 \begin{array}{l}
 \displaystyle \i\f{\partial\psi}{\partial t}=\hn{d}\psi\\[.4cm]
 \displaystyle \psi(0,\cdot)=\psi_0
 \end{array}
 \right.
\end{equation}
where $\hn{d}$ is now a nonlinear map with domain
\begin{multline}
\label{eq-domnon}
 \D(\hn{d}):=\left\{u\in\Ld{2}{d}:\exists\lambda>0,\,\exists q\in\C\:\text{ such that }\right.\\
 \left.u-\kappa_d(\beta|q|^{2\sigma})q\Gd{\lambda}{d}=:\phi_\lambda\in\Hd{2}{d}\:\text{ and }\:\phi_\lambda(0)=\td{\lambda}{d}{\beta|q|^{2\sigma}}q\right\}
\end{multline}
and action
\begin{equation}
 \label{eq-actnon}
 \hn{d}u=-\Delta\phi_\lambda-\kappa_d(\beta|q|^{2\sigma})q\lambda\Gd{\lambda}{d},\qquad\forall u\in\D(\hn{d}),
\end{equation}
with $\kappa_d$ still defined by \eqref{eq-theta} (an equivalent definition can be constructed from \eqref{eq-1d0}-\eqref{eq-2d0}-\eqref{eq-3d0}-\eqref{eq-actdue}). Consistently with the linear case, \emph{strong solutions} on an interval $[0,T]$ of \eqref{eq-domnon} have to be meant as functions in $C^0([0,T];\D(\hn{d}))\cap C^1([0,T];\Ld{2}{d})$ which satisfy the former as an equality in $C^0([0,T];\Ld{2}{d})$ and the latter as an equality in $\D(\hn{d})$. Note that $\D(\hn{d})$ is no more a linear space, but nevertheless a notion of convergence is still defined there by a natural analogous of the graph norm.

Furthermore, it is worth mentioning that, since in this context no general theory is available, the study of \eqref{eq-cau-st} relies on the study of the nonlinear versions of \eqref{eq-lin-duhamel} and \eqref{eq-lincharge}, that is
\begin{equation}
 \label{eq-duhamel}
 \psi(t,x)=(\Ud{d}(t)\psi_0)(x)+\i\int_0^t\Ud{d}(t-s,x)\kappa_d(\beta|q(s)|^{2\sigma})q(s)\ds,
\end{equation}
and
\begin{equation}
 \label{eq-charge}
 q(t)+\int_0^t\Kd{d}(t-s)c_d(\beta|q(s)|^{2\sigma})q(s)\ds=m_df_d(t).
\end{equation}
More precisely, all the methods are based on proving that \eqref{eq-charge}, also known as \emph{nonlinear charge equation}, has a unique solution with sufficiently ``nice'' features so that the function $\psi$ defined by \eqref{eq-duhamel} is a solution of \eqref{eq-cau-st} with the required regularity.

%%%%%%%%%%%%%%%%%%%%%%%%%%%%%%%%%%%%%%%%%%%%%%%%%%%%%%%%%%%%%%%%%%%%%%%%%%%%%%%%
%%%%%%%%%%%%%%%%%%%%%%%%%%%%%%%%%%%%%%%%%%%%%%%%%%%%%%%%%%%%%%%%%%%%%%%%%%%%%%%%
%%%%%%%%%%%%%%%%%%%%%%%%%%%%%%%%%%%%%%%%%%%%%%%%%%%%%%%%%%%%%%%%%%%%%%%%%%%%%%%%
%%%%%%%%%%%%%%%%%%%%%%%%%%%%%%%%%%%%%%%%%%%%%%%%%%%%%%%%%%%%%%%%%%%%%%%%%%%%%%%%
%%%%%%%%%%%%%%%%%%%%%%%%%%%%%%%%%%%%%%%%%%%%%%%%%%%%%%%%%%%%%%%%%%%%%%%%%%%%%%%%
%%%%%%%%%%%%%%%%%%%%%%%%%%%%%%%%%%%%%%%%%%%%%%%%%%%%%%%%%%%%%%%%%%%%%%%%%%%%%%%%
%%%%%%%%%%%%%%%%%%%%%%%%%%%%%%%%%%%%%%%%%%%%%%%%%%%%%%%%%%%%%%%%%%%%%%%%%%%%%%%%
%%%%%%%%%%%%%%%%%%%%%%%%%%%%%%%%%%%%%%%%%%%%%%%%%%%%%%%%%%%%%%%%%%%%%%%%%%%%%%%%
%%%%%%%%%%%%%%%%%%%%%%%%%%%%%%%%%%%%%%%%%%%%%%%%%%%%%%%%%%%%%%%%%%%%%%%%%%%%%%%%
%%%%%%%%%%%%%%%%%%%%%%%%%%%%%%%%%%%%%%%%%%%%%%%%%%%%%%%%%%%%%%%%%%%%%%%%%%%%%%%%
%%%%%%%%%%%%%%%%%%%%%%%%%%%%%%%%%%%%%%%%%%%%%%%%%%%%%%%%%%%%%%%%%%%%%%%%%%%%%%%%
%%%%%%%%%%%%%%%%%%%%%%%%%%%%%%%%%%%%%%%%%%%%%%%%%%%%%%%%%%%%%%%%%%%%%%%%%%%%%%%%

\section{Local well-posedness and conservation laws}
\label{sec-LWP}

In this section we present results concerning the local well-posedness of \eqref{eq-cau-st} and the associated conservation laws. The first results on local well-posedness of \eqref{eq-cau-st} that appeared in the literature concerned existence and uniqueness of \emph{weak solutions}, where in this context a weak solution of \eqref{eq-cau-st} on an interval $[0,T]$ is a function $\psi$ such that
\begin{gather}
 \psi(t,\cdot)=\phi_\lambda(t,\cdot)+q(t)\Gd{\lambda}{d}\in V_d, \qquad\forall t\in[0,T],\qquad\text{and}\\[.2cm]
\label{eq-cau-weak}
 \left\{
 \begin{array}{rl}
 \displaystyle \i\f{d}{dt}\langle\chi,\psi(t,\cdot)\rangle_{\Ld{2}{d}}= & \langle\nabla\chi_\lambda,\nabla\phi_\lambda(t,\cdot)\rangle_{\Ld{2}{d}}+\lambda(\langle\chi_\lambda,\phi_\lambda(t,\cdot)\rangle_{\Ld{2}{d}}-\langle\chi,\psi(t,\cdot)\rangle_{\Ld{2}{d}})\\[.4cm]
 & +\:\td{\lambda}{d}{\beta|q(t)|^{2\sigma}}q_\chi^*q(t),\qquad\forall \chi=\chi_\lambda+q_\chi\Gd{\lambda}{d}\in V_d,\quad\forall t\in[0,T],\\[.4cm]
 \displaystyle \psi(0,\cdot)=\psi_0,\quad\text{in}\quad V_d. & 
 \end{array}
 \right.
\end{gather}
Given this definition one can state the first well-posedness result.

\begin{theorem}[$d=1$ in \cite{AT-01}, $d=2$ in \cite{CCT-19}, $d=3$ in \cite{ADFT-03}]
\label{thm-LWPweak}
 Let $\beta\in\R\setminus\{0\}$, $\sigma>0$ and $\psi_0\in V_d$. Therefore:
 \begin{enumerate}[label=(\roman*)]
  \item if $d=1,\,3$, then there exists $T>0$ such that \eqref{eq-cau-st} admits a unique weak solution on $[0,T]$;
  \item if $d=2$, $\sigma\geq1/2$ and $(1+|k|^\ep)\F(\psi_0)\in\Ld{1}{2}$, for some $\ep>0$, then there exists $T>0$ such that \eqref{eq-cau-st} admits a unique weak solution on $[0,T]$
 \end{enumerate}
 In addition, letting 
 \begin{equation}
  \label{eq-existencetime}
  \tn{\beta}{\sigma}{d}(\psi_0):=\sup\{T>0:\text{there exists a unique weak solution of \eqref{eq-cau-st} on }\:[0,T]\},
 \end{equation}
 there results that
 \begin{equation}
  \label{eq-massenergy}
  M(t)=M(0)\qquad\text{and}\qquad \en{d}(t)=\en{d}(0),\qquad\forall t\in[0,\tn{\beta}{\sigma}{d}(\psi_0)),
 \end{equation}
 where $M(t)$ is the \emph{mass} associated with $\psi(t,\cdot)$, defined by \eqref{eq-mass}, and $\en{d}(t)$ is the \emph{energy} associated with $\psi(t,\cdot)$, defined by
 \begin{multline}
  \label{eq-energy}
  \en{d}(t)=\en{d}(\psi(t,\cdot)):=\\[.4cm]\left\{
  \begin{array}{ll}
   \displaystyle \frac12\|\partial_x\psi(t,\cdot)\|_{\Lu{2}}^2+\frac{\beta|q(t)|^{2\sigma+2}}{2\sigma+2}, & \text{if}\quad d=1\\[.4cm]
   \displaystyle \frac12\|\nabla\phi_\lambda(t,\cdot)\|_{\Ld{2}{d}}^2+\f{\lambda}{2}(\|\phi_\lambda(t,\cdot)\|_{\Ld{2}{d}}^2-\|\psi(t,\cdot)\|_{\Ld{2}{d}}^2)+\f{\td{\lambda}{d}{\f{\beta|q(t)|^{2\sigma}}{\sigma+1}}|q(t)|^2}{2}, & \text{if}\quad d=2,\,3.
  \end{array}
  \right.
 \end{multline}
\end{theorem}

\begin{remark}
 Conservations of mass and energy are not surprising since the invariances pointed out in the linear case are still present in the nonlinear case. Moreover, one can easily see that, setting $\beta=\alpha$ and $\sigma=0$, there results $\en{d}(t)=\frac{1}{2}\qd{d}(\psi(t,\cdot))$. As a consequence, in $d=3$ it is possible to find a form of the energy where the parameter $\lambda$ does not appear.
\end{remark}

\begin{remark}
 In fact, as a byproduct of the proof of Theorem \ref{thm-LWPweak}, one has that the weak solution of \eqref{eq-cau-st} belongs to $C^0([0,T];V_d)$, for any $T\in(0,\tn{\beta}{\sigma}{d}(\psi_0))$.
\end{remark}

The further assumptions required in the two-dimensional case are mainly due to the different qualitative behavior of the integral kernel of the charge equation $\Kd{d}$, defined by \eqref{eq-kernel}. Indeed, whereas in odd dimensions the $\f12$-\emph{Abel kernel} possesses nice regularizing properties in Sobolev spaces, in the two-dimensional case the \emph{Volterra function of order -1} has no regularizing properties in those spaces due to its highly singular behavior at the origin (for details see \cite{CCT-19,CFT-17}). Moreover, while in odd dimension Theorem \ref{thm-LWPweak} can be extended to finitely many nonlinear delta potentials, in the two-dimensional case the possibility of such an extension is open (and far from being understood).

Finally, again in odd dimension, it is possible to strengthen the results of Theorem \ref{thm-LWPweak}. Note that the following statement does not coincide with the ones contained in the original sources. It is written in such a way to be consistent with the notation used in Theorem \ref{thm-LWPweak}.

\begin{theorem}[$d=1$ in {\cite[case $s=1$]{CFT-19}}, $d=3$ in \cite{CFNT-17}]
\label{thm-LWPstrong}
 Let $d=1,\,3$, $\beta\in\R\setminus\{0\}$, $\sigma>0$ and $\psi_0\in \D(\hn{d})$. Then, there exists a unique strong solution of \eqref{eq-cau-st} on $[0,T]$, for every $T\in(0,\tn{\beta}{\sigma}{d}(\psi_0))$ (with $\tn{\beta}{\sigma}{d}(\psi_0)$ defined by \eqref{eq-existencetime}).
\end{theorem}

\begin{remark}
 In $d=2$ such a result on strong solutions is still missing due again to the features of the kernel of the charge equation, mentioned above.
\end{remark}

%%%%%%%%%%%%%%%%%%%%%%%%%%%%%%%%%%%%%%%%%%%%%%%%%%%%%%%%%%%%%%%%%%%%%%%%%%%%%%%%
%%%%%%%%%%%%%%%%%%%%%%%%%%%%%%%%%%%%%%%%%%%%%%%%%%%%%%%%%%%%%%%%%%%%%%%%%%%%%%%%
%%%%%%%%%%%%%%%%%%%%%%%%%%%%%%%%%%%%%%%%%%%%%%%%%%%%%%%%%%%%%%%%%%%%%%%%%%%%%%%%
%%%%%%%%%%%%%%%%%%%%%%%%%%%%%%%%%%%%%%%%%%%%%%%%%%%%%%%%%%%%%%%%%%%%%%%%%%%%%%%%
%%%%%%%%%%%%%%%%%%%%%%%%%%%%%%%%%%%%%%%%%%%%%%%%%%%%%%%%%%%%%%%%%%%%%%%%%%%%%%%%
%%%%%%%%%%%%%%%%%%%%%%%%%%%%%%%%%%%%%%%%%%%%%%%%%%%%%%%%%%%%%%%%%%%%%%%%%%%%%%%%
%%%%%%%%%%%%%%%%%%%%%%%%%%%%%%%%%%%%%%%%%%%%%%%%%%%%%%%%%%%%%%%%%%%%%%%%%%%%%%%%
%%%%%%%%%%%%%%%%%%%%%%%%%%%%%%%%%%%%%%%%%%%%%%%%%%%%%%%%%%%%%%%%%%%%%%%%%%%%%%%%
%%%%%%%%%%%%%%%%%%%%%%%%%%%%%%%%%%%%%%%%%%%%%%%%%%%%%%%%%%%%%%%%%%%%%%%%%%%%%%%%
%%%%%%%%%%%%%%%%%%%%%%%%%%%%%%%%%%%%%%%%%%%%%%%%%%%%%%%%%%%%%%%%%%%%%%%%%%%%%%%%
%%%%%%%%%%%%%%%%%%%%%%%%%%%%%%%%%%%%%%%%%%%%%%%%%%%%%%%%%%%%%%%%%%%%%%%%%%%%%%%%
%%%%%%%%%%%%%%%%%%%%%%%%%%%%%%%%%%%%%%%%%%%%%%%%%%%%%%%%%%%%%%%%%%%%%%%%%%%%%%%%

\section{Global well-posedness}
\label{sec-GWP}

The study of global well-posedness consists of detecting those cases in which the quantity $\tn{\beta}{\sigma}{d}(\psi_0)$ defined by \eqref{eq-existencetime} is equal to $+\infty$; namely, the cases in which the solution is \emph{global-in-time}.

In this context, the crucial points usually are:
\begin{enumerate}[label=(\roman*)]
 \item proving the so-called \emph{blow-up alternative}, i.e.
\begin{equation}
\label{eq-bua}
 \tn{\beta}{\sigma}{d}(\psi_0)<+\infty\qquad\text{and}\qquad \limsup_{t\nearrow \tn{\beta}{\sigma}{d}(\psi_0)}\|\psi(t,\cdot)\|_{V_d}<+\infty\qquad\Longrightarrow\qquad\text{contradiction}
\end{equation}
or equivalently
\begin{equation}
\label{eq-buaq}
 \qquad \tn{\beta}{\sigma}{d}(\psi_0)<+\infty\qquad\text{and}\qquad \limsup_{t\nearrow \tn{\beta}{\sigma}{d}(\psi_0)}|q(t)|<+\infty\qquad\Longrightarrow\qquad\text{contradiction};
\end{equation}
\item proving that conservations of mass and energy entail a-priori boundedness for $\|\psi(t,\cdot)\|_{V_d}$, or $|q(t)|$ on $[0,\tn{\beta}{\sigma}{d}(\psi_0))$.
\end{enumerate}

For the sake of clarity we distinguish the result on the \emph{defocusing} case, i.e. $\beta>0$, from those on the \emph{focusing} case, i.e. $\beta<0$.

\begin{theorem}[$d=1$ in \cite{AT-01}, $d=2$ in \cite{CCT-19}, $d=3$ in \cite{ADFT-03}]
 \label{thm-GWPd}
 Let $\beta,\,\sigma>0$ and $\psi_0\in V_d$ and let $\tn{\beta}{\sigma}{d}(\psi_0)$ be defined as in \eqref{eq-existencetime}. Therefore:
 \begin{enumerate}[label=(\roman*)]
  \item if $d=1,\,3$, then $\tn{\beta}{\sigma}{d}(\psi_0)=+\infty$;
  \item if $d=2$, $\sigma\geq1/2$ and $(1+|k|^\ep)\F(\psi_0)\in\Ld{1}{2}$, for some $\ep>0$, then $\tn{\beta}{\sigma}{2}(\psi_0)=+\infty$.
 \end{enumerate}
 Moreover, in these cases the weak solution of \eqref{eq-cau-st} belongs to $L^\infty([0,+\infty);V_d)$.
\end{theorem}

\begin{remark}
 Note that the presence of further assumptions when $d=2$ in Theorem \ref{thm-GWPd} is only due to the fact that they are necessary in Theorem \ref{thm-LWPweak}.  
\end{remark}

\begin{theorem}[$d=1$ in \cite{AT-01}, $d=3$ in \cite{ADFT-03}]
 \label{thm-GWPf}
 Let $d=1,3$, $\beta<0$, $\sigma\in(0,1]$ and $\psi_0\in V_d$ and let $\tn{\beta}{\sigma}{d}(\psi_0)$ be defined as in \eqref{eq-existencetime}. Therefore:
 \begin{enumerate}[label=(\roman*)]
  \item if $\sigma<1$, then $\tn{\beta}{\sigma}{d}(\psi_0)=+\infty$;
  \item if $\sigma=1$, then there exists $\mu_\beta^d>0$ such that $\tn{\beta}{1}{d}(\psi_0)=+\infty$ whenever $M(\psi_0)<\mu_\beta^d$.
 \end{enumerate}
 Moreover, in these cases the weak solution of \eqref{eq-cau-st} belongs to $L^\infty([0,+\infty);V_d)$.
\end{theorem}

The limitations on $\sigma$ displayed by Theorem \ref{thm-GWPf} are sharp, in the sense that, as we will see in Section \ref{sec-blow-up}, beyond the threshold $\sigma=1$, solutions which are not global-in-time may arise. As a consequence, the power $\sigma=1$ is usually called $L^2$-\emph{critical} power, while smaller powers are called $L^2$-\emph{subcritical} powers and larger powers are called $L^2$-\emph{supercritical} powers. In addition, also the threshold on the mass $\mu_\beta^d$, present in the case $\sigma=1$, is sharp in the same sense as before. As a consequence, such parameter is called $L^2$-\emph{critical mass}. It can be in fact explicitly computed and reads:
\begin{equation}
 \label{eq-crmass}
 \mu_\beta^d:=\left\{
 \begin{array}{ll}
  \displaystyle -\f{2}{\beta}, & \text{if}\quad d=1,\\[.4cm]
  \displaystyle \f{1}{-32\pi^2\beta}, & \text{if}\quad d=3.
 \end{array}
 \right.
\end{equation}

Finally, it is worth mentioning that also the absence of $d=2$ in Theorem \ref{thm-GWPf} is sharp. Indeed, as we will see again in Section \ref{sec-blow-up}, in the two-dimensional case solutions which are not global-in-time may be found for any value of $\sigma$. In other words, the reason for which there is no $L^2$-critical power in $d=2$ is that in this case any power is $L^2$-supercritical. Consistently, the concept of $L^2$-critical mass is not defined in this case.

\begin{remark}
 Concerning the extension of the previous theorems to the case of finitely many nonlinear deltas, this is known for Theorem \ref{thm-GWPd}(i), for Theorem \ref{thm-GWPf}(i) and, in the case $d=1$, for Theorem \ref{thm-GWPf}(ii).
\end{remark}

%%%%%%%%%%%%%%%%%%%%%%%%%%%%%%%%%%%%%%%%%%%%%%%%%%%%%%%%%%%%%%%%%%%%%%%%%%%%%%%%
%%%%%%%%%%%%%%%%%%%%%%%%%%%%%%%%%%%%%%%%%%%%%%%%%%%%%%%%%%%%%%%%%%%%%%%%%%%%%%%%
%%%%%%%%%%%%%%%%%%%%%%%%%%%%%%%%%%%%%%%%%%%%%%%%%%%%%%%%%%%%%%%%%%%%%%%%%%%%%%%%
%%%%%%%%%%%%%%%%%%%%%%%%%%%%%%%%%%%%%%%%%%%%%%%%%%%%%%%%%%%%%%%%%%%%%%%%%%%%%%%%
%%%%%%%%%%%%%%%%%%%%%%%%%%%%%%%%%%%%%%%%%%%%%%%%%%%%%%%%%%%%%%%%%%%%%%%%%%%%%%%%
%%%%%%%%%%%%%%%%%%%%%%%%%%%%%%%%%%%%%%%%%%%%%%%%%%%%%%%%%%%%%%%%%%%%%%%%%%%%%%%%
%%%%%%%%%%%%%%%%%%%%%%%%%%%%%%%%%%%%%%%%%%%%%%%%%%%%%%%%%%%%%%%%%%%%%%%%%%%%%%%%
%%%%%%%%%%%%%%%%%%%%%%%%%%%%%%%%%%%%%%%%%%%%%%%%%%%%%%%%%%%%%%%%%%%%%%%%%%%%%%%%
%%%%%%%%%%%%%%%%%%%%%%%%%%%%%%%%%%%%%%%%%%%%%%%%%%%%%%%%%%%%%%%%%%%%%%%%%%%%%%%%
%%%%%%%%%%%%%%%%%%%%%%%%%%%%%%%%%%%%%%%%%%%%%%%%%%%%%%%%%%%%%%%%%%%%%%%%%%%%%%%%
%%%%%%%%%%%%%%%%%%%%%%%%%%%%%%%%%%%%%%%%%%%%%%%%%%%%%%%%%%%%%%%%%%%%%%%%%%%%%%%%
%%%%%%%%%%%%%%%%%%%%%%%%%%%%%%%%%%%%%%%%%%%%%%%%%%%%%%%%%%%%%%%%%%%%%%%%%%%%%%%%

\section{Standing waves}
\label{sec-standing}

A further point of interest in the study of \eqref{eq-cau-st} is given by the \emph{standing waves}, i.e. solutions of the equation of the form $\psin{\beta}{\sigma}{d}(t,x)=\e^{\i\omega t}\un{\beta}{\sigma}{d}(x)$, for some $\omega\in\R$ usually called \emph{frequency}. Clearly, the search for such solutions reduces to the search for the functions
\[
 \un{\beta}{\sigma}{d}\in\D(\hn{d})\qquad\text{that satisfy}\qquad(\hn{d}+\omega)\un{\beta}{\sigma}{d}=0.
\]
These are usually called \emph{bound states} of the NLS with concentrated nonlinearity, can be completely classified in any dimension and are positive up to the \emph{gauge invariance} $\e^{\i\theta}$, $\theta\in[0,2\pi)$.

\begin{theorem}[$d=1$ in \cite{BD-21} and {\cite[case $s=1$]{CFT-19}}, $d=2$ in \cite{ACCT-21}, $d=3$ in \cite{ANO-13}]
 \label{thm-bound}
 Let $\sigma>0$. Up to gauge invariance, the bound states of the NLS with concentrated nonlinearity are of the form
 \begin{equation}
  \label{eq-bound}
  \un{\beta}{\sigma}{d}(x):=q_{\beta,\sigma}^d(\omega)\Gd{\omega}{d}(x)
 \end{equation}
 with
 \begin{equation}
  \label{eq-parameter}
  q_{\beta,\sigma}^d(\omega):=\left\{
  \begin{array}{ll}
   \displaystyle \left(\f{2^{2\sigma+1}\omega^{\sigma+1/2}}{-\beta}\right)^{1/2\sigma}, & \text{if}\quad d=1,\\[.4cm]
   \displaystyle \left(\f{\log(\sqrt{\omega}/2)+\gamma}{-2\pi\beta}\right)^{1/2\sigma}, & \text{if}\quad d=2,\\[.4cm]
   \displaystyle \left(\f{\sqrt{\omega}}{-4\pi\beta}\right)^{1/2\sigma}, & \text{if}\quad d=3,
  \end{array}
  \right.
 \end{equation}
 where
 \begin{itemize}
  \item if $d=1,\,3$, then $\omega$ varies in $(0,+\infty)$ and $\beta$ varies in $(-\infty,0)$; while,
  \item if $d=2$, then $\omega$ varies in $(0,4\e^{-2\gamma})$ whenever $\beta$ varies in $(0,+\infty)$, and in $(4\e^{-2\gamma},+\infty)$ whenever $\beta$ varies in $(-\infty,0)$.
 \end{itemize}
\end{theorem}

Note that the difference between the odd and the even dimensions is remarkable. Indeed, in $d=2$ there is a branch of bound states also in the defocusing case $\beta>0$. This is a phenomenon which possesses an analogous neither in the context of the NLS with concentrated nonlinearity nor in the context of the standard NLS.

\begin{remark}
\label{rem-qualitative}
 In fact, it is also possible to compute explicitly $\en{d}(\un{\beta}{\sigma}{d})$ and $M(\un{\beta}{\sigma}{d})$ as functions of the frequency $\omega$. As a consequence, one can establish the following qualitative behaviors (that are relevant in the study of the orbital stability discussed below). If $d=1,\,3$, then
  \begin{itemize}
   \item whenever $\sigma<1$
   \begin{itemize}
   \item $\en{d}(\un{\beta}{\sigma}{d})$ is negative, continuous and decreasing on $\R^+$ and
   \[
    \lim_{\omega\searrow0}\en{d}(\un{\beta}{\sigma}{d})=0,\qquad \lim_{\omega\to+\infty}\en{d}(\un{\beta}{\sigma}{d})=-\infty,
   \]
   \item $M(\un{\beta}{\sigma}{d})$ is continuous and increasing on $\R^+$ and
   \[
    \lim_{\omega\searrow0}M(\un{\beta}{\sigma}{d})=0,\qquad \lim_{\omega\to+\infty}M(\un{\beta}{\sigma}{d})=+\infty;
   \]
   \end{itemize}
   \item whenever $\sigma>1$
   \begin{itemize}
   \item $\en{d}(\un{\beta}{\sigma}{d})$ is positive, continuous and increasing on $\R^+$ and
   \[
    \lim_{\omega\searrow0}\en{d}(\un{\beta}{\sigma}{d})=0,\qquad\lim_{\omega\to+\infty}\en{d}(\un{\beta}{\sigma}{d})=+\infty,
   \]
   \item $M(\un{\beta}{\sigma}{d})$ is continuous and decreasing on $\R^+$ and
   \[
    \lim_{\omega\searrow0}M(\un{\beta}{\sigma}{d})=+\infty, \qquad \lim_{\omega\to+\infty}M(\un{\beta}{\sigma}{d})=0;
   \]
   \end{itemize}
   \item  whenever $\sigma=1$, $E^d_{\beta,1}(\un{\beta}{1}{d})=0$ and $M(\un{\beta}{1}{d})=\mu_\beta^d$ (with $\mu_\beta^d$ defined by \eqref{eq-crmass}), for all $\omega>0$.
  \end{itemize}
   If on the contrary $d=2$, then
  \begin{itemize}
  \item whenever $\beta<0$
  \begin{itemize}
  \item $\en{2}(\un{\beta}{\sigma}{2})$ is continuous on $(4\e^{-2\gamma},+\infty)$, decreasing on $(4\e^{-2\gamma},4\e^{-2\gamma+1/\sigma}]$ and increasing on $[4\e^{-2\gamma+1/\sigma},+\infty)$, and
  \[
   \lim_{\omega\searrow4\e^{-2\gamma}}\en{2}(\un{\beta}{\sigma}{2})=0, \qquad \lim_{\omega\to+\infty}\en{2}(\un{\beta}{\sigma}{2})=+\infty,
  \]
  \item $M(\un{\beta}{\sigma}{2})$ is continuous on $(4\e^{-2\gamma},+\infty)$, increasing on $(4\e^{-2\gamma},4\e^{-2\gamma+1/\sigma}]$ and decreasing on $[4\e^{-2\gamma+1/\sigma},+\infty)$, and
  \[
   \lim_{\omega\searrow4\e^{-2\gamma}}M(\un{\beta}{\sigma}{2})=\lim_{\omega\to+\infty}M(\un{\beta}{\sigma}{2})=0;
  \]
  \end{itemize}
  \item whenever $\beta>0$
  \begin{itemize}
  \item $\en{2}(\un{\beta}{\sigma}{2})$ is continuous and increasing on $(0,4\e^{-2\gamma})$, and
  \[
    \lim_{\omega\searrow0}\en{2}(\un{\beta}{\sigma}{2})=-\infty, \qquad \lim_{\omega\nearrow4\e^{-2\gamma}}\en{2}(\un{\beta}{\sigma}{2})=0,
  \]
  \item $M(\un{\beta}{\sigma}{2})$ is continuous and decreasing on $(0,4\e^{-2\gamma})$, and
  \[
   \lim_{\omega\searrow0}M(\un{\beta}{\sigma}{2})=+\infty, \qquad \lim_{\omega\nearrow4\e^{-2\gamma}}M(\un{\beta}{\sigma}{2})=0.
  \]
  \end{itemize}
  \end{itemize}
\end{remark}

%%%%%%%%%%%%%%%%%%%%%%%%%%%%%%%%%%%%%%%%%%%%%%%%%%%%%%%%%%%%%%%%%%%%%%%%%%%%%%%%
%%%%%%%%%%%%%%%%%%%%%%%%%%%%%%%%%%%%%%%%%%%%%%%%%%%%%%%%%%%%%%%%%%%%%%%%%%%%%%%%
%%%%%%%%%%%%%%%%%%%%%%%%%%%%%%%%%%%%%%%%%%%%%%%%%%%%%%%%%%%%%%%%%%%%%%%%%%%%%%%%
%%%%%%%%%%%%%%%%%%%%%%%%%%%%%%%%%%%%%%%%%%%%%%%%%%%%%%%%%%%%%%%%%%%%%%%%%%%%%%%%
%%%%%%%%%%%%%%%%%%%%%%%%%%%%%%%%%%%%%%%%%%%%%%%%%%%%%%%%%%%%%%%%%%%%%%%%%%%%%%%%
%%%%%%%%%%%%%%%%%%%%%%%%%%%%%%%%%%%%%%%%%%%%%%%%%%%%%%%%%%%%%%%%%%%%%%%%%%%%%%%%

\subsection{Orbital stability}
\label{subsec-orbital}

The most natural question concerning the bound states is about their \emph{stability}. In particular, one starts discussing \emph{orbital} stability.

Precisely, a bound state $\un{\beta}{\sigma}{d}$ is said \emph{orbitally stable} when, for every $\ep>0$, there exists $\delta>0$ such that: if $\|\psi_0-\e^{i\theta}\un{\beta}{\sigma}{d}\|_{V_d}<\delta$ for some $\theta\in[0,2\pi)$ and $\psi$ is a weak solution of \eqref{eq-cau-st} on $[0,\tn{\beta}{\sigma}{d}(\psi_0))$, then $\psi$ can be continued to a solution on $[0,+\infty)$ and
\[
 \sup_{t\in\R^+}\inf_{\theta\in[0,2\pi)}\|\psi(t,\cdot)-\e^{i\theta}\un{\beta}{\sigma}{d}\|_{V_d}<\ep.
\]
Otherwise $\un{\beta}{\sigma}{d}$ is said \emph{orbitally unstable}.

The results on orbital stability established in this context are (mainly) obtained by using the methods introduced by \cite{GSS-87} (and the qualitative behaviors pointed our in Remark \ref{rem-qualitative}). Preliminarily, set
\begin{equation}
 \label{eq-aux}
 a_d:=\left\{
 \begin{array}{ll}
  \displaystyle 0, & \text{if}\quad d=1,\,3,\\[.4cm]
  \displaystyle 4\e^{-2\gamma}, & \text{if}\quad d=2,
 \end{array}
 \right.
\end{equation}

\begin{theorem}[$d=1$ in \cite{BD-21}, $d=2$ in \cite{ACCT-21}, $d=3$ in \cite{ANO-13}]
 \label{thm-ostab}
 Let $\sigma>0$ and let $\beta<0$, when $d=1,\,3$, and $\beta\in\R\setminus\{0\}$, when $d=2$. Therefore,
 \begin{enumerate}[label=(\roman*)]
  \item if $d=1,3$ and $\sigma<1$, then $\un{\beta}{\sigma}{d}$ is orbitally stable for every $\omega>a_d$; while, if $d=2$, $\beta>0$ and $\sigma\geq1/2$, then $\un{\beta}{\sigma}{2}$ is orbitally stable for every $\omega\in(0,a_2)$;
  \item if $d=3$ and $\sigma\geq1$ or if $d=2$, $\beta<0$ and $\sigma\geq1/2$, then $\un{\beta}{\sigma}{d}$ is orbitally unstable for every $\omega>a_d$.
 \end{enumerate}
 \end{theorem}
 
 First note that the case $d=1,\,\sigma\geq1$ is not explicitly addressed by the literature. However, this lack is not significant in the sense that it is not due to challenges in the proof, but rather to the fact that historically the discussion on the one-dimensional case mainly focused on other issues. In this case, the result and proof should be analogous to those of the dimension three.
 
 Note also that, again, the features of the dimension two are remarkable. Here, the transition between stability and instability does not occur at some specific power, i.e. the $L^2$-critical power, but switching from defocusing to focusing. As a consequence, here any focusing power nonlinearity can be legitimately considered $L^2$-supercritical. On the contrary, the threshold $\sigma\geq1/2$ is not significant since, again, it is connected to local well-posedness and is purely technical.
 
 Finally, in strict connection with orbital stability, one can wonder in which cases the bound states are actually \emph{ground states} at some fixed mass $\mu>0$ of the NLS with concentrated nonlinearity, i.e. functions $u\in V_d$ such that $\|u\|_{\Ld{2}{d}}^2=\mu$ and
 \[
  \en{d}(u)=\En{d}(\mu):=\inf_{\substack{v\in V_d\\\|v\|_{\Ld{2}{d}}^2=\mu}}\en{d}(v).
 \]
 
 \begin{theorem}[$d=1$ in \cite{BD-21}, $d=2$ in \cite{ACCT-21}]
 \label{thm-gstate}
 Let $\sigma>0$ and let $\beta<0$, when $d=1$, and $\beta\in\R\setminus\{0\}$, when $d=2$. Therefore,
 \begin{enumerate}[label=(\roman*)]
  \item if $d=1$ and $\sigma<1$ or if $d=2$ and $\beta>0$, then $\En{d}(\mu)\in(-\infty,0)$, for all $\mu>0$, and $\un{\beta}{\sigma}{d}$ is the unique (up to gauge invariance) ground state at mass $M(\un{\beta}{\sigma}{d})$ of the NLS with concentrated nonlinearity, for every $\omega>a_1$ in the former case and every $\omega\in(0,a_2)$ in the latter case;
  \item if $d=1$ and $\sigma=1$, then 
  \[
   \mathcal{E}_{\beta,1}^1(\mu)=\left\{
   \begin{array}{ll}
    0, & \text{if}\quad\mu\in(0,\mu_\beta^1],\\[.2cm]
    -\infty, & \text{if}\quad\mu>\mu_\beta^1,
   \end{array}
   \right.
  \]
  and all the bound states $(\un{\beta}{1}{1})_{\omega>0}$ are actually ground states at mass $\mu_\beta^1$ of the NLS with concentrated nonlinearity;
  \item if $d=1$ and $\sigma>1$ or if $d=2$ and $\beta<0$, then $\En{d}(\mu)=-\infty$, for all $\mu>0$.
 \end{enumerate}
 \end{theorem}

 Combining the results of Theorems \ref{thm-bound} and \ref{thm-gstate} and of Remark \ref{rem-qualitative} one sees that:
 \begin{itemize}
  \item in case (i), there exists a unique positive ground state at mass $\mu$, for every $\mu>0$;
  \item in case (ii), positive ground states at mass $\mu$ exist if and only if $\mu=\mu_\beta^1$ and are infinitely many;
  \item in case (iii) no bound state is in fact a ground state.
 \end{itemize}
 Finally, we mention that the previous result has never been explicitly discussed for $d=3$. However, it should be possible to prove a result analogous to the case $d=1$, using the same techniques.

%%%%%%%%%%%%%%%%%%%%%%%%%%%%%%%%%%%%%%%%%%%%%%%%%%%%%%%%%%%%%%%%%%%%%%%%%%%%%%%%
%%%%%%%%%%%%%%%%%%%%%%%%%%%%%%%%%%%%%%%%%%%%%%%%%%%%%%%%%%%%%%%%%%%%%%%%%%%%%%%%
%%%%%%%%%%%%%%%%%%%%%%%%%%%%%%%%%%%%%%%%%%%%%%%%%%%%%%%%%%%%%%%%%%%%%%%%%%%%%%%%
%%%%%%%%%%%%%%%%%%%%%%%%%%%%%%%%%%%%%%%%%%%%%%%%%%%%%%%%%%%%%%%%%%%%%%%%%%%%%%%%
%%%%%%%%%%%%%%%%%%%%%%%%%%%%%%%%%%%%%%%%%%%%%%%%%%%%%%%%%%%%%%%%%%%%%%%%%%%%%%%%
%%%%%%%%%%%%%%%%%%%%%%%%%%%%%%%%%%%%%%%%%%%%%%%%%%%%%%%%%%%%%%%%%%%%%%%%%%%%%%%%

\subsection{Asymptotic stability}
\label{subsec-asymptotic}

Another perspective on stability is that of the \emph{asymptotic stability}, which roughly speaking means that, if the initial datum of \eqref{eq-cau-st} is ``close'' to a bound state, then the solution gets asymptotically ``closer and closer'' to a, possibly distinct, bound state.

The results mentioned below only deal with the case $d=3$. Analogous results are obtained for $d=1$ in \cite{BKKS-08,KKS-12}, but for another class of nonlinearities, which is very general from a certain point of view, but does not take into account powers. Nevertheless, it is reasonable to guess that one should obtain the same outcomes also for powers. Nothing is known, instead, for $d=2$.

In order to state the results, some preliminary definitions are required. First, we set $L^1_w(\R^3):=L^1(\R^3,(1+|x|^{-1})\dx)$. Moreover, let $L_{\omega}^\sigma:\Ld{2}{3}\times\Ld{2}{3}\to\Ld{2}{3}\times\Ld{2}{3}$ be the operator defined by
\[
 L_{\omega}^\sigma:=
 \begin{pmatrix}
  0 & \mathcal{H}_{\alpha_1(\sigma,\omega)}^3+\omega\\
  -\mathcal{H}_{\alpha_2(\sigma,\omega)}^3-\omega & 0
 \end{pmatrix}
\]
where
\[
 \alpha_1(\sigma,\omega):=-(2\sigma+1)\f{\sqrt{\omega}}{4\pi},\qquad\text{and}\qquad\alpha_2(\sigma,\omega):=-\f{\sqrt{\omega}}{4\pi},
\]
which is a byproduct of the linearization of $\hn{3}$ around a bound state of frequency $\omega$. Such an operator possesses two purely imaginary eigenvalues $\pm\i\xi_{\sigma,\omega}$, with $\xi_{\sigma,\omega}:=2\sigma\omega\sqrt{1-\sigma^2}$ (only defined for $\sigma\leq1$), whose eigenspaces are spanned by
\begin{equation}
\label{eq-aux2}
 \Phi_{\sigma,\omega}(x):=\Gd{{\omega-\xi_{\sigma,\omega}}}{3}(x)\begin{pmatrix} 1 \\ \i \end{pmatrix}-\f{\sqrt{1-\sigma^2}-1}{\sigma}\Gd{{\omega+\xi_{\sigma,\omega}}}{3}(x)\begin{pmatrix} 1 \\ -\i \end{pmatrix}
\end{equation}
and its complex conjugate. Let us also denote by $\Phi_{\sigma,\omega}^1$ and $\Phi_{\sigma,\omega}^2$ the first and the second components of $\Phi_{\sigma,\omega}$, respectively.

\begin{theorem}
 \label{thm-astab}
 Let $d=3$, $\beta<0$, $\omega>0$ and $\theta\in[0,2\pi)$. Therefore:
 \begin{enumerate}[label=(\roman*)]
  \item (\cite{ANO-13}) if $\sigma\in(0,1/\sqrt{2})$ and $\psi_0\in V_3\cap L^1_w(\R^3)$ with $\|\psi_0-\e^{\i\theta}\un{\beta}{\sigma}{3}\|_{V_3\cap L^1_w(\R^3)}$ sufficiently small, then there exist $\omega_\infty>0$, $\psi_\infty:\R^3\to\C$ and $r_\infty:\R^+\times\R^3\to\C$, with $\psi_\infty,\,r_\infty(t)\in\Ld{2}{3}$, for all $t>0$, and $\|r_\infty(t)\|_{\Ld{2}{3}}=O(t^{-5/4})$ as $t\to+\infty$, such that the weak solution $\psi$ of \eqref{eq-cau-st} satisfies
  \[
   \psi(t,x)=\e^{\i\omega_\infty t}u_{\beta,\sigma,\omega_{\infty}}^3(x)+(\Ud{3}(t)\psi_\infty)(x)+r_\infty(t,x),\qquad\text{as}\quad t\to+\infty;
  \]
  \item (\cite{ANO-16}) there exists $\displaystyle\sigma^*\in\bigg(\f{1}{\sqrt{2}},\f{\sqrt{3}+1}{2\sqrt{2}}\bigg]$ such that, if $\sigma\in(1/\sqrt{2},\sigma^*)$ and $\psi_0\in V_3\cap L^1_w(\R^3)$ and satisfies
  \[
   \psi_0=\e^{\i\omega+\vartheta}\un{\beta}{\sigma}{3}+\e^{\i\omega+\vartheta}\left[(z+\overline{z})\Phi_{\sigma,\omega}^1+\i(z-\overline{z})\Phi_{\sigma,\omega}^2\right]+f
  \]
  for some $\vartheta\in\R$, $z\in\C$ and $f\in\Ld{2}{3}\cap L^1_w(\R^3)$ such that
  \[
   |z|\leq\sqrt{\ep},\qquad\text{and}\qquad\|f\|_{L^1_w(\R^3)}\leq C\ep^{3/2}
  \]
  with $C>0$ and $\ep>0$ small, then there exist $\omega_\infty,\,k_\infty>0$, $b\in\R$, $\psi_\infty:\R^3\to\C$ and $r_\infty:\R^+\times\R^3\to\C$, with $\psi_\infty,\,r_\infty(t)\in\Ld{2}{3}$, for all $t>0$, and $\|r_\infty(t)\|_{\Ld{2}{3}}=O(t^{-1/4})$  as $t\to+\infty$, such that the weak solution $\psi$ of \eqref{eq-cau-st} satisfies
  \[
   \psi(t,x)=\e^{\i\omega_\infty t+\i b\log(1+\ep k_\infty t)}u_{\beta,\sigma,\omega_{\infty}}^3(x)+(\Ud{3}(t)\psi_\infty)(x)+r_\infty(t,x),\qquad\text{as}\quad t\to+\infty.
  \]
 \end{enumerate}
\end{theorem}

%%%%%%%%%%%%%%%%%%%%%%%%%%%%%%%%%%%%%%%%%%%%%%%%%%%%%%%%%%%%%%%%%%%%%%%%%%%%%%%%
%%%%%%%%%%%%%%%%%%%%%%%%%%%%%%%%%%%%%%%%%%%%%%%%%%%%%%%%%%%%%%%%%%%%%%%%%%%%%%%%
%%%%%%%%%%%%%%%%%%%%%%%%%%%%%%%%%%%%%%%%%%%%%%%%%%%%%%%%%%%%%%%%%%%%%%%%%%%%%%%%
%%%%%%%%%%%%%%%%%%%%%%%%%%%%%%%%%%%%%%%%%%%%%%%%%%%%%%%%%%%%%%%%%%%%%%%%%%%%%%%%
%%%%%%%%%%%%%%%%%%%%%%%%%%%%%%%%%%%%%%%%%%%%%%%%%%%%%%%%%%%%%%%%%%%%%%%%%%%%%%%%
%%%%%%%%%%%%%%%%%%%%%%%%%%%%%%%%%%%%%%%%%%%%%%%%%%%%%%%%%%%%%%%%%%%%%%%%%%%%%%%%
%%%%%%%%%%%%%%%%%%%%%%%%%%%%%%%%%%%%%%%%%%%%%%%%%%%%%%%%%%%%%%%%%%%%%%%%%%%%%%%%
%%%%%%%%%%%%%%%%%%%%%%%%%%%%%%%%%%%%%%%%%%%%%%%%%%%%%%%%%%%%%%%%%%%%%%%%%%%%%%%%
%%%%%%%%%%%%%%%%%%%%%%%%%%%%%%%%%%%%%%%%%%%%%%%%%%%%%%%%%%%%%%%%%%%%%%%%%%%%%%%%
%%%%%%%%%%%%%%%%%%%%%%%%%%%%%%%%%%%%%%%%%%%%%%%%%%%%%%%%%%%%%%%%%%%%%%%%%%%%%%%%
%%%%%%%%%%%%%%%%%%%%%%%%%%%%%%%%%%%%%%%%%%%%%%%%%%%%%%%%%%%%%%%%%%%%%%%%%%%%%%%%
%%%%%%%%%%%%%%%%%%%%%%%%%%%%%%%%%%%%%%%%%%%%%%%%%%%%%%%%%%%%%%%%%%%%%%%%%%%%%%%%

\section{Blow-up analysis}
\label{sec-blow-up}

The last question on \eqref{eq-cau-st} discussed by the literature is the so called \emph{blow-up analysis}, that is the search for solutions that \emph{blow-up} in a finite time. By \eqref{eq-existencetime}, \eqref{eq-bua} and \eqref{eq-buaq} this is equivalent to prove that $\tn{\beta}{\sigma}{d}(\psi_0)<+\infty$.

%%%%%%%%%%%%%%%%%%%%%%%%%%%%%%%%%%%%%%%%%%%%%%%%%%%%%%%%%%%%%%%%%%%%%%%%%%%%%%%%
%%%%%%%%%%%%%%%%%%%%%%%%%%%%%%%%%%%%%%%%%%%%%%%%%%%%%%%%%%%%%%%%%%%%%%%%%%%%%%%%
%%%%%%%%%%%%%%%%%%%%%%%%%%%%%%%%%%%%%%%%%%%%%%%%%%%%%%%%%%%%%%%%%%%%%%%%%%%%%%%%
%%%%%%%%%%%%%%%%%%%%%%%%%%%%%%%%%%%%%%%%%%%%%%%%%%%%%%%%%%%%%%%%%%%%%%%%%%%%%%%%
%%%%%%%%%%%%%%%%%%%%%%%%%%%%%%%%%%%%%%%%%%%%%%%%%%%%%%%%%%%%%%%%%%%%%%%%%%%%%%%%
%%%%%%%%%%%%%%%%%%%%%%%%%%%%%%%%%%%%%%%%%%%%%%%%%%%%%%%%%%%%%%%%%%%%%%%%%%%%%%%%

\subsection{First results}
\label{subsec-seminal}

First blow-up results are based on a method, introduced by \cite{G-77}, which essentially relies on the detection of the properties of the so called \emph{moment of inertia} associated with a solution $\psi$ of \eqref{eq-cau-st} with initial datum $\psi_0$, i.e.
\[
 \M{d}(t):=\int_{\R^d}|x\,\psi(t,x)|^2\dx.
\]
Precisely, in this context the main point is proving that
\begin{equation}
\label{eq-virial}
 \f{d^2\M{d}(t)}{dt^2}=8\en{d}(\psi_0)+g_{\beta,\sigma}^d(|q(t)|^2),
\end{equation}
with
\[
 g_{\beta,\sigma}^d(y):=\left\{
 \begin{array}{ll}
  \displaystyle 4\beta\f{\sigma-1}{\sigma+1}y^{\sigma+1}, & \text{if}\quad d=1,\,3,\\[.4cm]
  \displaystyle 2\bigg(\f{1}{\pi}-\f{4\beta\sigma y^\sigma}{\sigma+1}\bigg)\,y, & \text{if}\quad d=2.
 \end{array}
 \right.
\]
Equality \eqref{eq-virial} is usually called \emph{Virial Identity} and enables one to establish, in some cases, uniform concavity of $\M{d}$ depending on the energy of the initial datum. Since positivity and uniform concavity are not consistent with globality-in-time it is clear how this method allows to detect the cases in which blow-up solutions arise.

\begin{remark}
 Note that by Theorem \ref{thm-GWPd} blow-up analysis is meaningful only when $\beta<0$ and, for $d=1,\,3$, when $\sigma\geq1$.
\end{remark}

Preliminarily, recall that a \emph{Schwartz function} is a $C^\infty(\R^d)$-function which decays faster than the reciprocal of any polinomial as $|x|\to+\infty$ and whose derivatives of any order decay faster than the reciprocal of any polinomial as $|x|\to+\infty$.

\begin{theorem}[$d=1$ in \cite{AT-01}, $d=2$ in \cite{ACCT-20}, $d=3$ in \cite{ADFT-04}]
 \label{thm-bu}
 Let $\beta<0$, $\sigma>0$ and $\psi_0\in V_d$ and let $\tn{\beta}{\sigma}{d}(\psi_0)$ be defined as in \eqref{eq-existencetime}. Moreover, assume that:
 \begin{enumerate}[label=(\roman*)]
  \item when $d=1,\,3$, $\sigma\geq1$; while,
  \item when $d=2$, $\sigma\geq1/2$ and the regular part of $\psi_0$ is a Schwartz function.
 \end{enumerate}
 Then
 \begin{equation}
  \label{eq-bucond}
  \en{d}(\psi_0)<\inf_{\omega>a_d}\en{d}(\un{\beta}{\sigma}{d})\qquad\Longrightarrow\qquad \tn{\beta}{\sigma}{d}(\psi_0)<+\infty,
 \end{equation}
 where $(\un{\beta}{\sigma}{d})_{\omega>a_d}$ is the family of the bound states given by Theorem \ref{thm-bound} and $a_d$ is defined by \eqref{eq-aux}.
\end{theorem}

Some comments are in order. First, $\inf_{\omega>a_d}\en{d}(\un{\beta}{\sigma}{d})$ can be explicitly computed and there results
\[
 \inf_{\omega>a_d}\en{d}(\un{\beta}{\sigma}{d})=\left\{
 \begin{array}{ll}
  \displaystyle 0 & \text{if} \quad d=1,\,3,\\[.2cm]
  \displaystyle -\f{\sigma}{4\pi(\sigma+1)(-4\pi\sigma\beta)^{1/\sigma}} & \text{if} \quad d=2.
 \end{array}
 \right.
\]
On the other hand, \eqref{eq-bucond} is sharp in the sense that above this energy threshold one has, for instance, standing waves, which are global-in-time.

Moreover, we underline that in $d=2$ the assumptions on $\sigma$ and on the regularity of the initial datum are just technical. On the contrary, it is remarkable that in this case one can find blowing-up solutions for any focusing power nonlinearity (again, a phenomenon with no analogue neither in the context of the NLS with concentrated nonlinearity nor in the context of the standard NLS).

Finally, note that, when $d=1,\,3$ and $\sigma=1$, \eqref{eq-bucond} does not contradict point (ii) of Theorem \ref{thm-GWPf} since one can actually prove (using some suitable versions of the Gagliardo-Nirenberg inequalities established in \cite{AT-01} for $d=1$ and in \cite{ADFT-03} for $d=3$) that in these cases
\[
  E_{\beta,1}^d(\psi_0)<\inf_{\omega>0}E_{\beta,1}^d(\un{\beta}{1}{d})=0\qquad\Longrightarrow\qquad M(\psi_0)>\mu_\beta^d,
\]
where $\mu_\beta^d$ is again the $L^2$-critical mass defined by \eqref{eq-crmass}.

%%%%%%%%%%%%%%%%%%%%%%%%%%%%%%%%%%%%%%%%%%%%%%%%%%%%%%%%%%%%%%%%%%%%%%%%%%%%%%%%
%%%%%%%%%%%%%%%%%%%%%%%%%%%%%%%%%%%%%%%%%%%%%%%%%%%%%%%%%%%%%%%%%%%%%%%%%%%%%%%%
%%%%%%%%%%%%%%%%%%%%%%%%%%%%%%%%%%%%%%%%%%%%%%%%%%%%%%%%%%%%%%%%%%%%%%%%%%%%%%%%
%%%%%%%%%%%%%%%%%%%%%%%%%%%%%%%%%%%%%%%%%%%%%%%%%%%%%%%%%%%%%%%%%%%%%%%%%%%%%%%%
%%%%%%%%%%%%%%%%%%%%%%%%%%%%%%%%%%%%%%%%%%%%%%%%%%%%%%%%%%%%%%%%%%%%%%%%%%%%%%%%
%%%%%%%%%%%%%%%%%%%%%%%%%%%%%%%%%%%%%%%%%%%%%%%%%%%%%%%%%%%%%%%%%%%%%%%%%%%%%%%%

\subsection{Pseudoconformal invariance}
\label{subsec-pseudo}

In addition to Theorem \ref{thm-bu}, whenever $d=1,\,3$, $\beta<0$ and $\sigma=1$ it is possible to explicitly construct blowing-up solutions whose initial datum satisfies $M(\psi_0)=\mu_\beta^d$. This can be done thanks to an additional symmetry arising in the $L^2$-critical case, the so called \emph{pseudoconformal invariance}, which can be stated as follows.

\begin{proposition}[$d=1$ in \cite{AT-99}, $d=3$ in \cite{ADFT-04}]
\label{prop-inv}
 Let $d=1,3$, $\beta\in\R\setminus\{0\}$ and $\sigma=1$. If $\psi$ is a solution of \eqref{eq-cau-st} with initial datum $\psi_0$, then, for every fixed $T>0$,
 \[
  \Psi_T^d(t,x):=\f{\e^{-\i\f{|x|^2}{4(T-t)}}}{(T-t)^{d/2}}\psi\bigg(\f{1}{T-t},\f{x}{T-t}\bigg)
 \]
 is a solution of \eqref{eq-cau-st} with initial datum
 \[
  \Psi_{T,0}^d(x):=\f{\e^{-\i\f{|x|^2}{4T}}}{T^{d/2}}\psi\bigg(\f{1}{T},\f{x}{T}\bigg).
 \]
\end{proposition}

Hence, it is sufficient to apply the previous result with $\psi$ equal to a standing wave to get blow-up at the $L^2$-critical mass.

\begin{theorem}[$d=1$ in \cite{AT-99,AT-01}, $d=3$ in \cite{ADFT-04}]
\label{thm-inv}
 Let $d=1,3$, $\beta<0$ and $\sigma=1$. Then, for every fixed $\omega>0$, $\theta\in[0,2\pi)$ and $T>0$, the function
 \[
  \Psi_{T,\beta,\omega}^d(t,x):=\f{\e^{-\i\f{|x|^2}{4(T-t)}+\i\f{\omega}{T-t}+\i\theta}}{(T-t)^{d/2}}\un{\beta}{1}{d}\bigg(\f{x}{T-t}\bigg),
 \]
 with $\un{\beta}{1}{d}$ a bound state given by Theorem \ref{thm-bound}, is a solution of \eqref{eq-cau-st} with initial datum
 \[
  \Psi_{T,0,\beta,\omega}^d(x):=\f{\e^{-\i\f{|x|^2}{4T}+\i\f{\omega}{T}+\i\theta}}{T^{d/2}}\un{\beta}{1}{d}\bigg(\f{x}{T}\bigg)
 \]
 that blows up at $T$. In other words $\tn{\beta}{1}{d}(\Psi_{T,0,\beta,\omega}^d)=T<+\infty$.
\end{theorem}

%%%%%%%%%%%%%%%%%%%%%%%%%%%%%%%%%%%%%%%%%%%%%%%%%%%%%%%%%%%%%%%%%%%%%%%%%%%%%%%%
%%%%%%%%%%%%%%%%%%%%%%%%%%%%%%%%%%%%%%%%%%%%%%%%%%%%%%%%%%%%%%%%%%%%%%%%%%%%%%%%
%%%%%%%%%%%%%%%%%%%%%%%%%%%%%%%%%%%%%%%%%%%%%%%%%%%%%%%%%%%%%%%%%%%%%%%%%%%%%%%%
%%%%%%%%%%%%%%%%%%%%%%%%%%%%%%%%%%%%%%%%%%%%%%%%%%%%%%%%%%%%%%%%%%%%%%%%%%%%%%%%
%%%%%%%%%%%%%%%%%%%%%%%%%%%%%%%%%%%%%%%%%%%%%%%%%%%%%%%%%%%%%%%%%%%%%%%%%%%%%%%%
%%%%%%%%%%%%%%%%%%%%%%%%%%%%%%%%%%%%%%%%%%%%%%%%%%%%%%%%%%%%%%%%%%%%%%%%%%%%%%%%

\subsection{Further analysis}
\label{subsec-further}

Concluding the section, we mention some finer results concerning blow-up analysis. Note that such results are available only for $d=1$, while everything is open for $d=2,3$. Note also that the following results only concern the focusing case and precisely $\beta=-1$, for the sake of simplicity.

The former concerns an upper bound for the \emph{blow-up rate} of blow-up solutions.

\begin{theorem}[\cite{HL-20}]
 \label{thm-burate}
 Let $d=-\beta=1$, $\sigma\geq1$ and $\psi_0\in V_1$ and let $\tn{-1}{\sigma}{1}(\psi_0)$ be defined as in \eqref{eq-existencetime}. If $\tn{-1}{\sigma}{1}(\psi_0)<+\infty$, then there exists $C_\sigma>0$ such that
 \[
  \|\partial_x\psi(t,\cdot)\|_{\Lu{2}}\geq \sqrt{|E(\psi_0)|}\qquad\Longrightarrow\qquad\|\partial_x\psi(t,\cdot)\|_{\Lu{2}}\geq C_\sigma(\tn{-1}{\sigma}{1}(\psi_0)-t)^{-\f{1-\sigma_c}{2}}
 \]
 with
 \begin{equation}
  \label{eq-sigmac}
  \sigma_c:=\f{\sigma+1}{2\sigma}.
 \end{equation}
\end{theorem}

Furthermore, it is possible to establish a sharp dichotomy between global existence and blow-up both in the $L^2$-supercritical and in the $L^2$-critical cases.

\begin{theorem}[\cite{HL-20}]
 \label{thm-dichotomy}
 Let $d=-\beta=1$, $\sigma\geq1$ and $\psi_0\in V_1$ and let $\tn{-1}{\sigma}{1}(\psi_0)$ be defined as in \eqref{eq-existencetime}.
 \begin{enumerate}[label=(\roman*)]
  \item On the one hand, set $\sigma>1$, assume
  \[
   M(\psi_0)^{\f{1-\sigma_c}{\sigma_c}}E(\psi_0)<M(\unu{-1}{\sigma}{1})^{\f{1-\sigma_c}{\sigma_c}}E(\unu{-1}{\sigma}{1})
  \]
  and define
  \[
   \eta(t):=\f{\|\psi_0\|_{\Lu{2}}^{\f{1-\sigma_c}{\sigma_c}}\|\partial_x\psi(t,\cdot)\|_{\Lu{2}}}{\|\unu{-1}{\sigma}{1}\|_{\Lu{2}}^{\f{1-\sigma_c}{\sigma_c}}\|\partial_x\unu{-1}{\sigma}{1}\|_{\Lu{2}}}.
  \]
  Therefore:
  \begin{itemize}
   \item if $\eta(0)<1$, then $\tn{-1}{\sigma}{1}(\psi_0)=+\infty$;
   \item if $\eta(0)>1$, then $\tn{-1}{\sigma}{1}(\psi_0)<+\infty$ and $\eta(t)>1$, for every $t\in[0,\tn{-1}{\sigma}{1}(\psi_0))$.
  \end{itemize}
  \item On the other hand, set $\sigma=1$. Therefore:
  \begin{itemize}
   \item if $M(\psi_0)<M(\unu{-1}{1}{1})=2$, then $E(\psi_0)>0$ and 
   \[
    \|\partial_x\psi(t,\cdot)\|_{\Lu{2}}\leq\f{2E(\psi_0)}{1-M(\psi_0)/2},\qquad\forall t\in[0,\tn{-1}{\sigma}{1}(\psi_0)),
   \]
    so that $\tn{-1}{\sigma}{1}(\psi_0)=+\infty$;
   \item if $E(\psi_0)<0$, then $\tn{-1}{\sigma}{1}(\psi_0)<+\infty$.
  \end{itemize}
 \end{enumerate}
\end{theorem}

\begin{remark}
 In fact, \cite{HL-20} presents also other results. It establishes some further regularity properties of the weak solutions of \eqref{eq-cau-st} and discusses both the phenomenon of mass concentration and the features of blow-up solutions at minimal and near-to-minimal mass in the $L^2$-critical case.
\end{remark}

It is also possible to establish a \emph{self-similar} structure for blow-up solutions in the $L^2$-supercritical case.

\begin{theorem}[\cite{HL-21}]
 \label{thm-selfsimilar}
 Let $d=-\beta=1$, $\sigma>1$ and let $\tn{-1}{\sigma}{1}(\psi_0)$ be defined as in \eqref{eq-existencetime}. Then, the function
 \[
  \psi(t,x):=\nu(t)^{1/2\sigma}\e^{\i\tau(t)}\zeta(\nu(t)x)
 \]
 weakly solves \eqref{eq-cau-st} (for some suitable $\psi_0\in V_1$) with $\lim_{t\nearrow T}\nu(t)=+\infty$, for some $T>0$, if and only if there exists $h>0$ and $k\in\R$ such that
 \[
  \nu(t)=\f{1}{\sqrt{2h(T-t)}},\qquad\tau(t)=\f{k}{2h}\log\bigg(\f{T}{T-t}\bigg)+\tau(0)
 \]
 and $\zeta(z)$ solves
 \[
  (k+\i h\sigma_c)\zeta-\i h\Lambda_z\zeta-\partial^2_{zz}\zeta-\delta_0|\zeta|^{2\sigma}\eta=0
 \]
 with $\sigma_c$ defined by \eqref{eq-sigmac}.
\end{theorem}

\begin{remark}
 Note that \cite{HL-21} also establishes some qualitative properties of the function $\zeta$.
\end{remark}

Finally, it is worth mentioning a result, which is more connected to \emph{scattering theory} in the cases where global existence is guaranteed, but which is nevertheless strictly related to the blow-up analysis. It concerns the so called \emph{asymptotic completeness}.

\begin{theorem}[\cite{AFH-21}]
 \label{thm-scattering}
 Let $d=-\beta=1$, $\sigma>1$ and $\psi_0\in V_1$. If
 \[
   M(\psi_0)^{\f{1-\sigma_c}{\sigma_c}}E(\psi_0)<M(\unu{-1}{\sigma}{1})^{\f{1-\sigma_c}{\sigma_c}}E(\unu{-1}{\sigma}{1}),
  \]
  and 
  \[
   \|\psi_0\|_{\Lu{2}}^{\f{1-\sigma_c}{\sigma_c}}\|\partial_x\psi_0\|_{\Lu{2}}<\|\unu{-1}{\sigma}{1}\|_{\Lu{2}}^{\f{1-\sigma_c}{\sigma_c}}\|\partial_x\unu{-1}{\sigma}{1}\|_{\Lu{2}}
  \]
  and $\psi$ is the weak solution of \eqref{eq-cau-st} with initial datum $\psi_0$, then there exists $\psi^+\in V_1$ such that
  \[
   \lim_{t\to+\infty}\|\Ud{1}(-t)\psi(t)-\psi^+\|_{V_1}=0.
  \]
\end{theorem}

\begin{remark}
 Note that a suitable analogue of Theorem \ref{thm-scattering} holds also in the defocusing case.
\end{remark}

%%%%%%%%%%%%%%%%%%%%%%%%%%%%%%%%%%%%%%%%%%%%%%%%%%%%%%%%%%%%%%%%%%%%%%%%%%%%%%%%
%%%%%%%%%%%%%%%%%%%%%%%%%%%%%%%%%%%%%%%%%%%%%%%%%%%%%%%%%%%%%%%%%%%%%%%%%%%%%%%%
%%%%%%%%%%%%%%%%%%%%%%%%%%%%%%%%%%%%%%%%%%%%%%%%%%%%%%%%%%%%%%%%%%%%%%%%%%%%%%%%
%%%%%%%%%%%%%%%%%%%%%%%%%%%%%%%%%%%%%%%%%%%%%%%%%%%%%%%%%%%%%%%%%%%%%%%%%%%%%%%%
%%%%%%%%%%%%%%%%%%%%%%%%%%%%%%%%%%%%%%%%%%%%%%%%%%%%%%%%%%%%%%%%%%%%%%%%%%%%%%%%
%%%%%%%%%%%%%%%%%%%%%%%%%%%%%%%%%%%%%%%%%%%%%%%%%%%%%%%%%%%%%%%%%%%%%%%%%%%%%%%%
%%%%%%%%%%%%%%%%%%%%%%%%%%%%%%%%%%%%%%%%%%%%%%%%%%%%%%%%%%%%%%%%%%%%%%%%%%%%%%%%
%%%%%%%%%%%%%%%%%%%%%%%%%%%%%%%%%%%%%%%%%%%%%%%%%%%%%%%%%%%%%%%%%%%%%%%%%%%%%%%%
%%%%%%%%%%%%%%%%%%%%%%%%%%%%%%%%%%%%%%%%%%%%%%%%%%%%%%%%%%%%%%%%%%%%%%%%%%%%%%%%
%%%%%%%%%%%%%%%%%%%%%%%%%%%%%%%%%%%%%%%%%%%%%%%%%%%%%%%%%%%%%%%%%%%%%%%%%%%%%%%%
%%%%%%%%%%%%%%%%%%%%%%%%%%%%%%%%%%%%%%%%%%%%%%%%%%%%%%%%%%%%%%%%%%%%%%%%%%%%%%%%
%%%%%%%%%%%%%%%%%%%%%%%%%%%%%%%%%%%%%%%%%%%%%%%%%%%%%%%%%%%%%%%%%%%%%%%%%%%%%%%%

\section{Derivation from the standard NLS equation}
\label{sec-derivation}

As mentioned in the introduction, a rigorous mathematical derivation of the NLS with concentrated nonlinearity from quantum many body systems in presence of impurities is still missing. What is available, at least for $d=1$ and, partially, for $d=3$, is a derivation from the standard NLS equation with suitably shrinking potentials. For $d=2$ nothing is known at the moment.

%%%%%%%%%%%%%%%%%%%%%%%%%%%%%%%%%%%%%%%%%%%%%%%%%%%%%%%%%%%%%%%%%%%%%%%%%%%%%%%%
%%%%%%%%%%%%%%%%%%%%%%%%%%%%%%%%%%%%%%%%%%%%%%%%%%%%%%%%%%%%%%%%%%%%%%%%%%%%%%%%
%%%%%%%%%%%%%%%%%%%%%%%%%%%%%%%%%%%%%%%%%%%%%%%%%%%%%%%%%%%%%%%%%%%%%%%%%%%%%%%%
%%%%%%%%%%%%%%%%%%%%%%%%%%%%%%%%%%%%%%%%%%%%%%%%%%%%%%%%%%%%%%%%%%%%%%%%%%%%%%%%
%%%%%%%%%%%%%%%%%%%%%%%%%%%%%%%%%%%%%%%%%%%%%%%%%%%%%%%%%%%%%%%%%%%%%%%%%%%%%%%%
%%%%%%%%%%%%%%%%%%%%%%%%%%%%%%%%%%%%%%%%%%%%%%%%%%%%%%%%%%%%%%%%%%%%%%%%%%%%%%%%

\subsection{The one-dimensional case}
\label{subsec-1d}

In this case the result that can be proven is exactly the expected one: the solution of \eqref{eq-cau-st} is the limit of the solution of a standard NLS with a suitably shrinking local potential.

Precisely, consider the following equation
\begin{equation}
 \label{eq-NLSapp}
 \i\f{\partial\psi_\ep}{\partial t}=-\f{\partial^2\psi_\ep}{\partial x^2}+\f{1}{\ep}V\bigg(\f{x}{\ep}\bigg)|\psi_\ep|^{2\sigma}\psi_\ep\qquad \sigma>0,\quad\ep>0,\qquad\text{on}\quad\R^+\times\R.
\end{equation}
It can be proved (see \cite{C-03}) that for any $\psi_0\in\Hu{1}$ and any $V\in L^1(\R,(1+|x|\dx))\cap \Lu{\infty}$,
\begin{itemize}
 \item if $V\geq0$ and $\sigma>0$, or
 \item if $V<0$ in at least an open interval and $\sigma\in(0,2)$ (or $\sigma=2$ and $\|\psi_0\|_{\Hu{1}}$ is small),
\end{itemize}
then the Cauchy problem associated with \eqref{eq-NLSapp} with initial datum $\psi_0$ admits a unique global \emph{strong} $H^1$\emph{-solution}, namely (see \cite{T-06}) a distributional solution in $C_{loc}^0([0,+\infty);\Hu{1})$ of
\begin{equation}
\label{eq-NLSappd}
 \psi_\ep(t,x)=(\Ud{1}(t)\psi_0)(x)-\f{\i}{\ep}\int_0^t\int_{\R}\Ud{1}(t-s,x-y)V\bigg(\f{y}{\ep}\bigg)|\psi_\ep(s,y)|^{2\sigma}\psi_\ep(s,y)\dy\ds,
\end{equation}
which is the \emph{Duhamel's formulation} of the Cauchy problem associated with \eqref{eq-NLSapp}. Then, we can state the following theorem.

\begin{theorem}[\cite{CFNT-14}]
 \label{thm-app1d}
 Let $\beta\in\R\setminus\{0\}$ and $V\in L^1(\R,(1+|x|\dx))\cap \Lu{\infty}$ with $\beta=\int_\R V(x)\dx$. Assume also that $V\geq0$ or $\sigma\in(0,1)$. Therefore, for every $\psi_0\in\Hu{1}$, if $\psi$ is the weak solution of \eqref{eq-cau-st} and $\psi_\ep$ is the strong $H^1$-solution of \eqref{eq-NLSappd}, then, for any fixed $T>0$,
 \[
  \sup_{t\in[0,T]}\|\psi_\ep(t,\cdot)-\psi(t,\cdot)\|_{\Hu{1}}\longrightarrow0,\qquad\text{as}\quad\ep\searrow0.
 \]
\end{theorem}

\begin{remark}
 In fact, in \cite{CFNT-14} the approximation established by Theorem \ref{thm-app1d} is proved in the case of finitely many nonlinear delta potentials by means of a standard NLS with finitely many shrinking local potentials.
\end{remark}

%%%%%%%%%%%%%%%%%%%%%%%%%%%%%%%%%%%%%%%%%%%%%%%%%%%%%%%%%%%%%%%%%%%%%%%%%%%%%%%%
%%%%%%%%%%%%%%%%%%%%%%%%%%%%%%%%%%%%%%%%%%%%%%%%%%%%%%%%%%%%%%%%%%%%%%%%%%%%%%%%
%%%%%%%%%%%%%%%%%%%%%%%%%%%%%%%%%%%%%%%%%%%%%%%%%%%%%%%%%%%%%%%%%%%%%%%%%%%%%%%%
%%%%%%%%%%%%%%%%%%%%%%%%%%%%%%%%%%%%%%%%%%%%%%%%%%%%%%%%%%%%%%%%%%%%%%%%%%%%%%%%
%%%%%%%%%%%%%%%%%%%%%%%%%%%%%%%%%%%%%%%%%%%%%%%%%%%%%%%%%%%%%%%%%%%%%%%%%%%%%%%%
%%%%%%%%%%%%%%%%%%%%%%%%%%%%%%%%%%%%%%%%%%%%%%%%%%%%%%%%%%%%%%%%%%%%%%%%%%%%%%%%

\subsection{The three-dimensional case}
\label{subsec-3d}

Unfortunately, in the three-dimensional case it is still open whether an analogous of Theorem \ref{thm-app1d} can be proved using local potentials. Nevertheless, it is possible to prove a version of it with a family of \emph{nonlocal} shrinking potentials.

Precisely, for every fixed $\sigma,\,\ep>0$ and $\beta\in\R\setminus\{0\}$, consider the following equation
\begin{equation}
 \label{eq-NLSappaux}
 \i\f{\partial\psi_\ep}{\partial t}=-\Delta\psi_\ep+\bigg(\f{-\ep}{\ell}+\beta\f{\ep^{2\sigma+2}}{\ell^{2\sigma+2}}|\langle\rho_\ep,\psi_\ep(t,\cdot)\rangle_{\Ld{2}{3}}|^{2\sigma}\bigg)\langle\rho_\ep,\psi_\ep(t,\cdot)\rangle_{\Ld{2}{3}}\rho_\ep,\qquad\text{on}\quad\R^+\times\R^3,
\end{equation}
where $(\rho_\ep)_{\ep>0}$ is defined by
\[
 \rho_\ep(x):=\f{1}{\ep^3}\rho\bigg(\f{x}{\ep}\bigg),
\]
with $\rho$ a real, positive, spherically symmetric Schwartz function of $\R^3$, (i.e. a kernel for $\delta_0$ in $\R^3$) and
\[
 \ell=\ell(\rho):=\langle\rho,(-\Delta)^{-1}\rho\rangle_{\Ld{2}{3}}.
\]
In this case, the \emph{Duhamel's formulation} of the Cauchy problem associated with \eqref{eq-NLSappaux} is
\begin{multline}
\label{eq-NLSappauxd}
 \psi_\ep(t,x)=(\Ud{3}(t)\psi\ep(0,\cdot))(x)+\\[.2cm]
 -\f{\i\ep}{\ell}\int_0^t\bigg(-1+\beta\f{\ep^{2\sigma+1}}{\ell^{2\sigma+1}}|\langle\rho_\ep,\psi_\ep(t,\cdot)\rangle_{\Ld{2}{3}}|^{2\sigma}\bigg)\langle\rho_\ep,\psi_\ep(t,\cdot)\rangle_{\Ld{2}{3}}(\Ud{3}(t-s)\rho_\ep)(x)\ds.
\end{multline}
One can prove (see \cite{CFNT-17}) that, whenever $\psi(0,\cdot)\in\Hd{2}{3}$, there exists a unique global $H^2$\emph{-strong solution} of \eqref{eq-NLSappaux}, namely a solution of \eqref{eq-NLSappauxd} in $C_{loc}^0([0,+\infty);\Hd{2}{3})\cap C_{loc}^1([0,+\infty);\Ld{2}{3})$.

Now, recall that in order to represent $\D(\hn{3})$ one can use the nonlinear analogous to \eqref{eq-3d0}, i.e.
\[
 \D(\hn{3})=\left\{u\in\Ld{2}{3}:\exists q\in\C\:\text{ s.t. }u-q\Gd{0}{3}=:\phi\in\Hdloc{2}{3}\cap\Ho{1}{3}\cap\Ho{2}{3},\:\phi(0)=\beta|q|^{2\sigma}q\right\}.
\]
Moreover, one can check that (\cite{CFNT-17})
 \begin{equation}
 \label{eq-appdat}
  \psi_{\ep,0}:=\phi_0+q_0\rho_\ep *\Gd{0}{3}\in\Hd{2}{3},\qquad\forall\psi_0\in\D(\hn{3}).
 \end{equation}
 
\begin{theorem}[\cite{CFNT-17}]
 \label{thm-app3d}
 Let $\beta\in\R\setminus\{0\}$ and $\ep>0$. Assume also that $\sigma>0$ when $\beta>0$, whereas $\sigma\in(0,1)$ when $\beta<0$. Therefore, for every $\psi_0\in\D(\hn{3})$, if $\psi$ is the weak solution of \eqref{eq-cau-st} and $\psi_\ep$ is the $H^2$-strong solution of \eqref{eq-NLSappaux} with $\psi_\ep(0,\cdot)=\psi_{\ep,0}$ ($\psi_{\ep,0}$ being defined by \eqref{eq-appdat}), then, for any fixed $T>0$,
 \[
  \sup_{t\in[0,T]}\|\psi_\ep(t,\cdot)-\psi(t,\cdot)\|_{\Ld{2}{3}}\leq C\ep^\delta,
 \]
 for some constant $C>0$ and $\delta\in(0,1/4)$.
\end{theorem}

%%%%%%%%%%%%%%%%%%%%%%%%%%%%%%%%%%%%%%%%%%%%%%%%%%%%%%%%%%%%%%%%%%%%%%%%%%%%%%%%
%%%%%%%%%%%%%%%%%%%%%%%%%%%%%%%%%%%%%%%%%%%%%%%%%%%%%%%%%%%%%%%%%%%%%%%%%%%%%%%%
%%%%%%%%%%%%%%%%%%%%%%%%%%%%%%%%%%%%%%%%%%%%%%%%%%%%%%%%%%%%%%%%%%%%%%%%%%%%%%%%
%%%%%%%%%%%%%%%%%%%%%%%%%%%%%%%%%%%%%%%%%%%%%%%%%%%%%%%%%%%%%%%%%%%%%%%%%%%%%%%%
%%%%%%%%%%%%%%%%%%%%%%%%%%%%%%%%%%%%%%%%%%%%%%%%%%%%%%%%%%%%%%%%%%%%%%%%%%%%%%%%
%%%%%%%%%%%%%%%%%%%%%%%%%%%%%%%%%%%%%%%%%%%%%%%%%%%%%%%%%%%%%%%%%%%%%%%%%%%%%%%%
%%%%%%%%%%%%%%%%%%%%%%%%%%%%%%%%%%%%%%%%%%%%%%%%%%%%%%%%%%%%%%%%%%%%%%%%%%%%%%%%
%%%%%%%%%%%%%%%%%%%%%%%%%%%%%%%%%%%%%%%%%%%%%%%%%%%%%%%%%%%%%%%%%%%%%%%%%%%%%%%%
%%%%%%%%%%%%%%%%%%%%%%%%%%%%%%%%%%%%%%%%%%%%%%%%%%%%%%%%%%%%%%%%%%%%%%%%%%%%%%%%
%%%%%%%%%%%%%%%%%%%%%%%%%%%%%%%%%%%%%%%%%%%%%%%%%%%%%%%%%%%%%%%%%%%%%%%%%%%%%%%%
%%%%%%%%%%%%%%%%%%%%%%%%%%%%%%%%%%%%%%%%%%%%%%%%%%%%%%%%%%%%%%%%%%%%%%%%%%%%%%%%
%%%%%%%%%%%%%%%%%%%%%%%%%%%%%%%%%%%%%%%%%%%%%%%%%%%%%%%%%%%%%%%%%%%%%%%%%%%%%%%%

\section{Results on connected problems}
\label{sec-connected}

Finally, it is worth reporting some works that address problems which are closely related in some sense to the NLS equation with point-concentrated nonlinearities.

%%%%%%%%%%%%%%%%%%%%%%%%%%%%%%%%%%%%%%%%%%%%%%%%%%%%%%%%%%%%%%%%%%%%%%%%%%%%%%%%
%%%%%%%%%%%%%%%%%%%%%%%%%%%%%%%%%%%%%%%%%%%%%%%%%%%%%%%%%%%%%%%%%%%%%%%%%%%%%%%%
%%%%%%%%%%%%%%%%%%%%%%%%%%%%%%%%%%%%%%%%%%%%%%%%%%%%%%%%%%%%%%%%%%%%%%%%%%%%%%%%
%%%%%%%%%%%%%%%%%%%%%%%%%%%%%%%%%%%%%%%%%%%%%%%%%%%%%%%%%%%%%%%%%%%%%%%%%%%%%%%%
%%%%%%%%%%%%%%%%%%%%%%%%%%%%%%%%%%%%%%%%%%%%%%%%%%%%%%%%%%%%%%%%%%%%%%%%%%%%%%%%
%%%%%%%%%%%%%%%%%%%%%%%%%%%%%%%%%%%%%%%%%%%%%%%%%%%%%%%%%%%%%%%%%%%%%%%%%%%%%%%%

\subsection{Nonlinear perturbations of the fractional Laplacian}
\label{subsec-fractional}

First, we mention a ``fractional perturbation'' in dimension one of what we introduced in Section \ref{sec-model} (nothing has been done in higher dimensions). The unique difference here is that the operator perturbed by a nonlinear Delta potential is the \emph{fractional Laplacian} $(-\Delta)^s$.

The problem is again defined by a suitable modification of the linear model (studied in \cite{MOS-18} and partially in \cite{CFT-19}), and reads exactly as \eqref{eq-cau-st} with $\hn{1}$ replaced by $\hs{1}$, where $\hs{1}$ is an operator whose definition is completely analogous to \eqref{eq-domnon}-\eqref{eq-actnon} with $\Gd{\lambda}{1}$ replaced by the Green's function of $(-\f{d^2}{dx^2})^s+\lambda$, denoted by $\Gs{\lambda}{1}$, and $-\Delta$ is replaced by $(-\f{d^2}{dx^2})^s$.

\begin{remark}
 As in the non fractional case, an analogous definition can be constructed from the fractional analogous to \eqref{eq-1d0}-\eqref{eq-2d0}-\eqref{eq-3d0}-\eqref{eq-actdue}.
\end{remark}

The problem has been addressed by \cite{CFT-19}, which establishes:
\begin{enumerate}[label=(\roman*)]
 \item local well-posedness in strong sense,
 \item conservation of mass and energy,
 \item global well-posedness in the defocusing and the focusing subcritical case (with critical power here depending on $s$);
 \item existence of blow-up solutions in the focusing supercritical case;
 \item complete classification of the standing waves.
\end{enumerate}

%%%%%%%%%%%%%%%%%%%%%%%%%%%%%%%%%%%%%%%%%%%%%%%%%%%%%%%%%%%%%%%%%%%%%%%%%%%%%%%%
%%%%%%%%%%%%%%%%%%%%%%%%%%%%%%%%%%%%%%%%%%%%%%%%%%%%%%%%%%%%%%%%%%%%%%%%%%%%%%%%
%%%%%%%%%%%%%%%%%%%%%%%%%%%%%%%%%%%%%%%%%%%%%%%%%%%%%%%%%%%%%%%%%%%%%%%%%%%%%%%%
%%%%%%%%%%%%%%%%%%%%%%%%%%%%%%%%%%%%%%%%%%%%%%%%%%%%%%%%%%%%%%%%%%%%%%%%%%%%%%%%
%%%%%%%%%%%%%%%%%%%%%%%%%%%%%%%%%%%%%%%%%%%%%%%%%%%%%%%%%%%%%%%%%%%%%%%%%%%%%%%%
%%%%%%%%%%%%%%%%%%%%%%%%%%%%%%%%%%%%%%%%%%%%%%%%%%%%%%%%%%%%%%%%%%%%%%%%%%%%%%%%

\subsection{Nonlinearity concentrated on a sphere}
\label{subsec-sfere}

Another variant of the problem introduced by Section \ref{sec-model} arises as one considers a non-zero-dimensional support for the delta measure. For instance, one may think to a smooth, compact manifold without boundary, such as, in particular, a sphere $\S\subset\R^3$.

Here, the whole setting is more technical, even in the linear case (see \cite{AGS-87} for $\S$ and \cite{BEL-14,BEHL-17,BLL-14} for more complex geometries). However, the proper variant of \eqref{eq-cau-st} is obtained by replacing $\hn{3}$ with the map $\hS{\S}$, with domain
\begin{align}
\label{eq-non_dom1}
\D (\hS{\S}): & = \big\{ u\in \Ld{2}{3}: \, \exists \lambda>0,\,q:\S\to\C \: \text{ s.t. } \\[.2cm]
\label{eq-non_dom2}           &  \hspace{5cm}u+\GPd{\lambda}{3}(\beta|q|^{2\sigma}q)=:\phi_\lambda \in \Hd{2}{3}\:\text{ and }\:q=u_{|\S}\big\}
\end{align}
and action
\begin{equation}
\label{eq-nonact}
\hS{\S} u =-\Delta\phi_\lambda+\lambda\GPd{\lambda}{3}(\beta|q|^{2\sigma}q),\qquad\forall u\in \D(\hS{\S}),
\end{equation}
where $\GPd{\lambda}{3}$ is now the \emph{Green's potential} associated with the unit sphere of the operator $-\Delta+\lambda$ in $\R^3$, i.e.
\begin{equation}
\label{eq-pot}
 \GPd{\lambda}{3}(h)[x]:=\int_{\S}\Gd{\lambda}{3}(x-y)\,  h(y)\,\dsy,\qquad \forall\, h:\S\to\C,\quad\forall x\in\R^3.
\end{equation}

\begin{remark}
 It is possible to find a definition in which the dumb parameter $\lambda$ does not appear by suitably adapting the sphere versions of \eqref{eq-1d0} and \eqref{eq-actdue}. Note that the link with the point-concentrated case in $d=1$ is due to the codimension of the support of the delta, which is one in both cases.
\end{remark}

The problem has been addressed by \cite{FTT-22}, which establishes:
\begin{enumerate}[label=(\roman*)]
 \item local well-posedness in strong sense,
 \item conservation of mass and energy,
 \item global well-posedness in the defocusing case for non large exponents of the nonlinearity.
\end{enumerate}

%%%%%%%%%%%%%%%%%%%%%%%%%%%%%%%%%%%%%%%%%%%%%%%%%%%%%%%%%%%%%%%%%%%%%%%%%%%%%%%%
%%%%%%%%%%%%%%%%%%%%%%%%%%%%%%%%%%%%%%%%%%%%%%%%%%%%%%%%%%%%%%%%%%%%%%%%%%%%%%%%
%%%%%%%%%%%%%%%%%%%%%%%%%%%%%%%%%%%%%%%%%%%%%%%%%%%%%%%%%%%%%%%%%%%%%%%%%%%%%%%%
%%%%%%%%%%%%%%%%%%%%%%%%%%%%%%%%%%%%%%%%%%%%%%%%%%%%%%%%%%%%%%%%%%%%%%%%%%%%%%%%
%%%%%%%%%%%%%%%%%%%%%%%%%%%%%%%%%%%%%%%%%%%%%%%%%%%%%%%%%%%%%%%%%%%%%%%%%%%%%%%%
%%%%%%%%%%%%%%%%%%%%%%%%%%%%%%%%%%%%%%%%%%%%%%%%%%%%%%%%%%%%%%%%%%%%%%%%%%%%%%%%

\subsection{Non-autonomous point-concentrated delta potentials}
\label{subsec-ionization}

Finally, we mention a class of problems where the point delta potentials are linear, but non-autonomous. This class of models is historically mentioned together with nonlinear delta potentials since, even though many features are different, the main techniques used for the investigation are the same.

Such models arise as one sets $\alpha=\alpha(t)$ in \eqref{eq-dom}-\eqref{eq-act}, with $\alpha:\R\to\C$ a given function, and considers the associated variant of \eqref{eq-cau-lin-st}. The resulting Cauchy problem involves a family of operators $\big(\Ht{d}{t}\big)_{t\in\R}$, indexed by the time variable $t$, which are self-adjoint at any fixed $t\in\R$. However, since the operator is not the same at every time, standard theory based on Stone's theorem does not apply any more, thus requiring to rely on some specific techniques used in the context of nonlinear delta potentials.

The main reason of interest for this model descends from an effective description of the microscopic dynamics of a quantum particle interacting with bosonic scalar quantum fields in configurations where the fields are very intense and the average number of carriers is large, in the so called \emph{quasi-classical limit} (see \cite{CCFO-21,CFO-19}). A typical example is given by the description of a quantum particle coupled with both acoustic and optical phonons in a compound ionic crystal (see, e.g., \cite[Chapter 4]{K-04}).

Global well-posedness for this problem has been studied in \cite{HMN-10} for $d=1$, in \cite{BCT-22} for $d=2$ and in \cite{SY-83}. Anyway, in this context, another interesting issue is the study of the asymptotic
behavior of the \emph{survival probability} of the $L^2$-normalized bound state $\varphi_{\alpha(0)}^d$ associated with the sole eigenvalue of $\Ht{d}{0}$ (see after \eqref{eq-eigen}), i.e.
\[
 \text{survival probability}:=|\langle\varphi_{\alpha(0)}^d,\psi(t,\cdot)\rangle_{\Ld{2}{d}}|^2
\]
where $\psi$ is the evolution of the initial state $\psi_0=\varphi_{\alpha(0)}^d$. Such a quantity can be proved to vanish at infinity, which shows a \emph{complete ionization} phenomenon, and an estimate of its decay rate can be also provided (\cite{CCLR-01,CCL-18} for $d=1$, \cite{BCT-22} for $d=2$ and \cite{CDFM-05} for $d=3$).

Finally, we underline that also a different non-autonomous delta model has been managed in the literature: the so-called \emph{traveling deltas} (\cite{DFT-00}). In this case, it is no more the strength of the delta potential which depends on time, but the the point where the delta measure is based.

%%%%%%%%%%%%%%%%%%%%%%%%%%%%%%%%%%%%%%%%%%%%%%%%%%%%%%%%%%%%%%%%%%%%%%%%%%%%%%%%
%%%%%%%%%%%%%%%%%%%%%%%%%%%%%%%%%%%%%%%%%%%%%%%%%%%%%%%%%%%%%%%%%%%%%%%%%%%%%%%%
%%%%%%%%%%%%%%%%%%%%%%%%%%%%%%%%%%%%%%%%%%%%%%%%%%%%%%%%%%%%%%%%%%%%%%%%%%%%%%%%
%%%%%%%%%%%%%%%%%%%%%%%%%%%%%%%%%%%%%%%%%%%%%%%%%%%%%%%%%%%%%%%%%%%%%%%%%%%%%%%%
%%%%%%%%%%%%%%%%%%%%%%%%%%%%%%%%%%%%%%%%%%%%%%%%%%%%%%%%%%%%%%%%%%%%%%%%%%%%%%%%
%%%%%%%%%%%%%%%%%%%%%%%%%%%%%%%%%%%%%%%%%%%%%%%%%%%%%%%%%%%%%%%%%%%%%%%%%%%%%%%%
%%%%%%%%%%%%%%%%%%%%%%%%%%%%%%%%%%%%%%%%%%%%%%%%%%%%%%%%%%%%%%%%%%%%%%%%%%%%%%%%
%%%%%%%%%%%%%%%%%%%%%%%%%%%%%%%%%%%%%%%%%%%%%%%%%%%%%%%%%%%%%%%%%%%%%%%%%%%%%%%%
%%%%%%%%%%%%%%%%%%%%%%%%%%%%%%%%%%%%%%%%%%%%%%%%%%%%%%%%%%%%%%%%%%%%%%%%%%%%%%%%
%%%%%%%%%%%%%%%%%%%%%%%%%%%%%%%%%%%%%%%%%%%%%%%%%%%%%%%%%%%%%%%%%%%%%%%%%%%%%%%%
%%%%%%%%%%%%%%%%%%%%%%%%%%%%%%%%%%%%%%%%%%%%%%%%%%%%%%%%%%%%%%%%%%%%%%%%%%%%%%%%
%%%%%%%%%%%%%%%%%%%%%%%%%%%%%%%%%%%%%%%%%%%%%%%%%%%%%%%%%%%%%%%%%%%%%%%%%%%%%%%%

\end{document}